\documentclass{aims}
\usepackage{txfonts}
\usepackage{mathrsfs}
\usepackage{eufrak}
\def\typeofarticle{Type of article}
\def\currentvolume{5}
\def\currentissue{x}
\def\currentyear{2018}

\def\ppages{xxx--xxx}
\def\DOI{}
\def\Received{July 2017}
\def\Accepted{September 2017}
\def\Published{December 2017}
\newcommand{\re}[1]{\textcolor{black}{#1}}

\newcommand{\dt}[1]{\textcolor{black}{#1}}

\numberwithin{equation}{section}

\usepackage{cleveref}
\def\bu{\textbf{u}}    
\def\bc{\textbf{c}} 
\def\R{\mathbb{R}}

\def\R{\mathbb{R}}
\def\N{\mathbb{N}}

\def\bu{\textbf{u}}
\def\bc{\textbf{c}}

\def\G{\mathcal{G}}

\def\Q{\mathcal{Q}}

\def\H{\mathcal{H}}

\def\E{\mathcal{E}}

\def\S{\mathcal{S}}

\def\M{\mathcal{M}}

\def\Mu{\mathfrak{M}}

\def\A{\mathcal{A}}

\def\restrict#1{\raise-.5ex\hbox{\ensuremath|}_{#1}}

\newcommand{\eequ}{\end{equation}}

\newcommand{\bequ}{\begin{equation}}

\newcommand{\eequd}{\end{eqnarray*}}

\newcommand{\bequd}{\begin{eqnarray*}}

 \graphicspath{{Figures/}}                   
\begin{document}
\title{Reconstruction \re{of Mutation Laws} in Heterogeneous Tumours with Local and Nonlocal Dynamics}
\author{%
  Maher Alwuthaynani\affil{1},
  Raluca Eftimie\affil{2}
  and
  Dumitru Trucu\affil{1,}\corrauth
}

\shortauthors{the Author(s)}

\address{%
  \addr{\affilnum{1}}{Division of Mathematics, University of Dundee, Dundee DD1 4HN, Scotland, UK;}
  \addr{\affilnum{2}}{Laboratoire Math\'ematiques de Besan\c con, UMR-CNRS 6623, Universit\'e de Bourgogne Franche-Comt\'e, 16 Route de Gray, Besan\c con, 25000, France;}
  }
%
\corraddr{trucu@maths.dundee.ac.uk}
%
%
\begin{abstract}
\re{Cancer cell mutations occur when cells undergo multiple cell divisions, and these mutations can be spontaneous or environmentally-induced. The mechanisms that promote and sustain these mutations are still not fully understood.}
This study deals with the identification (or reconstruction) of the usually unknown cancer cell mutation law, \re{which lead to the transformation of a primary tumour cell population into a secondary, more aggressive cell population. We focus on local and nonlocal mathematical models for cell dynamics and movement, and identify these mutation laws}  from macroscopic tumour snapshot data collected at some later stage in the tumour evolution. In \re{a} local cancer invasion model, we first reconstruct the mutation \re{law} when \re{we assume that} the mutations depend \re{only on} the \re{surrounding} cancer cells \re{(i.e., the ECM plays no role in mutations)}. Second, \re{we assume that the} mutations depend on the ECM only, \re{and we reconstruct the mutation law in this case}. Third, \re{we reconstruct} the mutation when we assume \re{that} there is no prior knowledge about the mutations. \re{Finally, for the} nonlocal cancer invasion model, we reconstruct the mutation law that depends on the cancer cells and on the ECM.  \re{For these numerical reconstructions, our} approximations \re{are} based on \re{the} finite difference method combined with \re{the} finite elements method. As the inverse problem is ill-posed, we use \re{the} Tikhonov regularisation technique in order to regularise the solution. Stability of the solution is examined by adding additive noise into the measurements.
\end{abstract}
\maketitle
%
%
%
%
%
%
\section{Introduction}
The beginning of a primary solid tumour is \re{the result of} a single normal cell that \re{is} transformed as a result of mutations in certain key genes. \re{Cells can mutate spontaneously, or mutations can be environmentally induced. Mutations occur during cell division, and most of the time the immune system can recognise mutated cells and eliminate them. When the immune system fails to eliminate cells with mutations in genes that control cell proliferation, cells become cancerous. It is known that all cancer cell lines have at least one mutation, with most cancer cell lines having more than one mutation; e.g., in \cite{Ikediobi2006_Mutations24genesCancer} 137 oncogenic mutations were identified in 14 out of 24 known cancer genes in 60 human cancer cell lines. Moreover, cancer cells are genetically unstable and cells inside the solid tumours keep mutating leading to very heterogeneous tumour masses. }
\par\re{The mechanisms behind the mutation pressure are still not fully understood. Experimental studies have shown that some changes in the extracellular matrix (ECM) can correlate with sustained cell proliferative signalling and an increased risk of developing cancer~\cite{Pickup2014_ECMmodulatesHallmarksCancer}. Other studies have shown that culturing cells for long times in stiff hydrogels can lead to the subclonal selection of genomic aberrations in cells~\cite{Lopez-Carrasco2020_ECMstiffnessCellMutation}, thus suggesting that the ECM properties could impact the mutation status of cells in solid tumours. Other studies suggested that the maintenance of cells at high density in the absence of proliferation leads to an increase in mutagenesis following cell division~\cite{JacksonLoeb}. Therefore, there seem to be different mechanisms that can trigger and influence the mutation rate of cells.
}
\par
\begin{figure}[htp]
	\centering
	\includegraphics[width=0.75\textwidth]{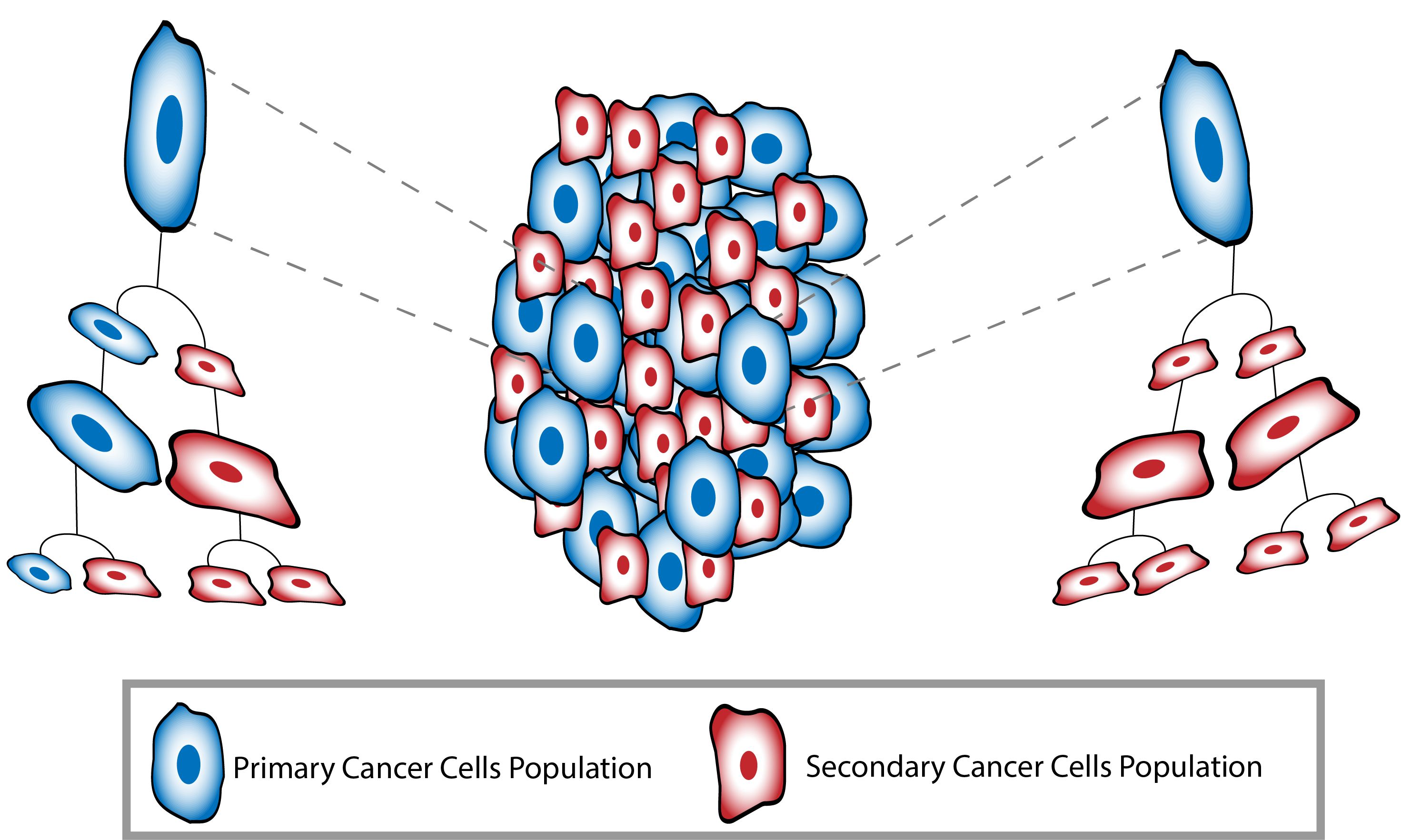}\\
	\caption{\small{\textit{Schematic of cancer cells proliferation.}}}
\label{fig:Schematic_One1}
\end{figure}

 \re{Mutated cells not only incur sustained proliferation (see Figure \ref{fig:Schematic_One1}), but can also exhibit migrational and invasion properties ~\cite{Novikov2021_MutationCancerMigrationInvasion}, eventually leading to cancer metastasis. The invasion of surrounding tissue is the result of ECM degradation and remodelling by the cancer cells (which can secrete various proteolytic enzymes, such as matrix metalloproteinases (MMPs)) as well as other cells in the microenvironment. The last decades have seen the development of numerous  theoretical studies based on mathematical models, which investigate computationally the biological mechanisms behind cancer cell invasion into the tissue. Most of these models are single-scale models; see~\cite{Anderson_et_al_2000, Anderson_2005,Anderson_Chaplain_1998,Gatenby_Gawlinski_1996,Chaplain_Lolas_2005,Domschke_et_al_2014} and references therein. More recently, multi-scale mathematical models have started to be developed to consider also the multi-scale aspects of various biological processes occurring during cancer invasion~\cite{Ramis-Conde_et_al_2008,Marciniak_Ptashnyk_2008,Macklin_et_al_2009,Deisboeck_et_al_2011,Dumitru_et_al_2013}. }
\par
\re{The main issue faced by all these single-scale and especially multi-scale mathematical models for cancer growth and invasion is parameter estimation. The last few years have seen the publication of various mathematical studies that try to estimate numerically different model parameters using inverse problem formulations~\cite{Colin2014_InverseproblemVascularisation,Gholami2016_InverseProblemGlioma,Hogea2008_InverseProblemTumGrowthMechanics,Jaroudi2018_InverseProblemLocationBrainTum,Subramanian2020_InverseProblemDiiffusionGrowth}. These studies focus on identifying the magnitude of isotropic~\cite{Subramanian2020_InverseProblemDiiffusionGrowth} and anisotropic tumour diffusion~\cite{Gholami2016_InverseProblemGlioma}, the magnitude of tumour growth rate \cite{Hogea2008_InverseProblemTumGrowthMechanics,Subramanian2020_InverseProblemDiiffusionGrowth}, the strength/location of tumour-induced tissue deformation~\cite{Hogea2008_InverseProblemTumGrowthMechanics}, the position of the blood vessels that act as a source for the oxygen concentration that influences tumour growth~\cite{Colin2014_InverseproblemVascularisation}, the location of the source of tumours~\cite{Jaroudi2018_InverseProblemLocationBrainTum}. }
\par
\re{To our knowledge, no studies have tried to estimate the mutation laws for cells inside heterogeneous tumours. To address this aspect, in this current study we consider a simplified problem with only two cancer cell populations (one mutating into the second one), that can exhibit both random and directed haptotactic movement. Moreover, we assume that the cells can degrade and remodel the surrounding ECM density. Furthermore, the type of interactions among cancer cells and between cells and ECM is not always very clear: some simpler mathematical models for cancer growth and invasion consider local cell-cell and cell-ECM interactions~\cite{Gatenby_Gawlinski_1996,Anderson_Chaplain_1998}, while other more complex models consider non-local cell-cell and cell-ECM interactions~\cite{Armstrong_et_al_2006,Dumitru_et_al_2013,Domschke_et_al_2014,Gerisch_Chaplain_2008}. In this study, we consider both type of models, local and non-local, and estimate the mutation laws for both cases. Finally, we assume here the knowledge of additional information in terms of both exact and noisy measurements of the tumour constituent density at some later time in the tumour evolution. We test our inversion approach on several cancer proliferation laws that are usually used in cancer modelling: logistic and Gompertz proliferation.
}

\re{This paper is structured as follows. In Section~\ref{A_mathematical_model_of_two_cancer_cell_sub_populations_I} we describe the local mathematical model for the dynamics of the two cancer cell sub-populations and the extracellular matrix (ECM). In Sections~\ref{Inverse_Problem_for_Unknown_Cancer_Cells_Mutation}-\ref{Inverse_Problem_Forward_Operator_and_Computational_Approach} we formulate the inverse problem for this local model under the assumptions of logistic cell proliferation and mutation depending only on cancer cell density. In Section~\ref{Recon_Mutation_in_Local_Model} we formulate the inverse problem for this local model under the assumption of logistic cell proliferation and mutation depending on both cell and ECM density. In Section~\ref{Sect6} we describe a non-local model for two cancer cell populations and their interactions with the ECM, while in Section~\ref{Numerical_Approach_for_Direct_Non_Local_Cancer_Model} we present the numerical approach for this forward nonlocal model. In Section~\ref{Sect8} we formulate the inverse problem for this nonlocal model for cancer invasion. We summarise and discuss our results in Section~\ref{Sect:Conclusion}. 
}

%
%
%
%
%
%
%
%
%
\section{Mathematical Model for Two Local Cancer Cell Sub-populations}
\label{A_mathematical_model_of_two_cancer_cell_sub_populations_I}
\re{In this study we consider two populations of tumour cells: a primary cell population $c_{1}(x,t)$ that can mutate into a more aggressive secondary cell population $c_{2}(x,t)$. These cell populations exercise a spatial redistribution via random movement (with diffusion coefficients $D_{1}$ and $D_{2}$) and directed movement towards extracellular matrix (ECM) gradients $v(x,t)$ (with haptotactic coefficients $\eta_{1}$ and $\eta_{2}$). Moreover, the two populations undergo logistic growth (at the same rate $\mu_{c}$), up to a carrying capacity ($K_{c}$).  We also assume that the ECM undergoes degradation (at a rate $\rho$) and remodelling (at a rate $\mu_{v}$). The above assumptions are described by the equations:}
%
\begin{subequations}
\label{eq:Direct_Mutation_DepOn_C}
	\begin{align}
	\label{eq:DP_C1_Simp_Mut_Chap4}
	\frac{\partial{c_{1}}} {\partial {t}} & = \underbrace{D_{1} {\Delta}  {c_{1}}}_{\text{diffusion}} - \underbrace{{\eta_{1}} {\nabla} \cdot \left( {c_{1}} {\nabla} {v} \right)}_{\text{haptotactic movement}} + \underbrace{\mu_{c} {c_{1}} \re{\left(1- \frac{c_{1} + c_{2} + v}{K_{c}} \right) } }_{\text{logistic proliferation}} - \dt{\underbrace{\omega(t)}_\text{\parbox[t]{0.4in}{mutation\\[-0.2cm] switch}}\underbrace{\Q(\cdot,\cdot)}_\text{\parbox[t]{0.4in}{unknown\\[-0.2cm] mutation}}},\\
	\label{eq:DP_C2_Simp_Mut_Chap4}
	\frac{\partial{c_{2}}} {\partial {t}} & = \underbrace{D_{2} {\Delta}  {c_{2}}}_{\text{diffusion}} - \underbrace{{\eta_{2}} {\nabla} \cdot \left( {c_{2}} {\nabla} {v} \right)}_{\text{haptotactic movement}} + \underbrace{\mu_{c} {c_{2}} \re{\left(1 - \frac{c_{1} + c_{2} + v}{K_{c}} \right) } }_{\text{logistic proliferation}} +  \dt{\underbrace{\omega(t)}_\text{\parbox[t]{0.4in}{mutation\\[-0.2cm] switch}}\underbrace{\Q(\cdot,\cdot)}_\text{\parbox[t]{0.4in}{unknown\\[-0.2cm] mutation}}},\\
	\label{eq:DP_V_Simp_Mut_Chap4}
	\frac{\partial {v}} {\partial {t} } & = \underbrace{- \rho \left(c_{1} + c_{2}\right) v}_{\text{degradation}} + \underbrace{\mu_{v} \left(K_{c} - v - c_{1} - c_{2} \right)^{\re{+}}}_{\text{remodelling term}}.
	\end{align}
\end{subequations}
\re{\dt{where} $(K_{c}-v-c_{1}-c_{2})^{+}:=max\{ (K_{c}-v-c_{1}-c_{2}),0\}$. Finally, the unknown term $\Q(\cdot,\cdot)$ \dt{represents} the mutation \dt{law} of cell subpopulation $c_{1}$ into cell subpopulation $c_{2}$, \dt{which is assumed to be mediated by a time-dependant mutation enhancement $\omega(t)$ that is known \emph{a priori} and is taken here of the form}
\bequd
\dt{\omega(t):= \frac{\bigg(1 + tanh\bigg(\frac{t-t_{1,2}}{t_{s}}\bigg)\bigg)}{2}},
\eequd
\dt{where} \re{$t_{1,2}$ is the time at which mutations from $c_{1}$ to $c_{2}$ start occurring, and $t_{s}>0$ is a time-steepness coefficient for this mutation law.} }\\

\re{\dt{The mutation law $\Q(\cdot, \cdot)$} is \dt{considerd here} unknown due to either unknown dependance on the primary cell population $c_{1}$, or unknown dependance on the ECM $v$, or unknown dependance on both primary tumour cell population and ECM}. \re{In this study we investigate three assumptions related to this mutation term, namely: (i) mutation depends linearly on the density of primary tumour cells; (ii) mutation depends linearly on the density of primary tumour cells, and nonlinearly on the ECM density~\cite{Lopez-Carrasco2020}; (iii) mutation law is very general and depends autonomously on the primary tumour and ECM. \dt{Thus mathematically, these cases correspond to three inverse problems that seek to identify the unknown mutation law $\Q(\cdot,\cdot)$ in the following three situations}:}
\begin{enumerate}
\item[(i)] \re{Mutation depends linearly on the density of primary tumour cell sub-population $c_{1}$ \dt{but does not depend at all on ECM}, and so this is given by the unknown term }
\bequ
\label{eq:Mutation_DepOn_C}
\dt{\Q(c_{1},v):=\widetilde{\Q}_{1}(c_{1})= \delta_{0} c_{1}},
\eequ
\dt{with} $\delta_{0}$ \dt{representing} the \dt{unknown} mutation rate.
 \item[(ii)] \re{Mutation depends \dt{in a known linear manner on $c_{1}$ and in an unknown nonlinear way on the density of ECM. The unknown dependence on $v$ is denoted mathematically by the unknown function $\widetilde{\Q}_{2}(v)$, and so the entire mutation law is therefore of the form}}
\begin{equation}
\label{eq:Mutation_DepOn_V}
\dt{\Q(c_{1},v) :=  \delta_{0} c_{1} \widetilde{Q}_{2}(v)},
\end{equation}
\dt{with $\delta_{0}>0$ here being considered known. A usual choice for $\widetilde{Q}_{2}(v)$ is of the form \cite{Shuttleworth_2020b,Shuttleworth_2020a}: }
\dt{\bequ\label{explicit_measurement_case_ii}
\widetilde{Q}_{2}(v):=
\left\{
\begin{array}{ll}
\frac{exp\left(\frac{-1}{\kappa^{2} - \left(1 - v\left(x, t\right)\right)^{2}}\right)}{exp\left(-\frac{1}{\kappa^{2}}\right)} & \quad \quad \text{if} \quad  1-\kappa < v\left(x, t\right) < 1,\\
0, &\quad \quad \text{otherwise}, 
\end{array}
\right.
\eequ}
\re{where $\kappa>0$ is a certain level of ECM beyond which mutations can occur.} 
\item[(iii)] \re{Mutation is given by an \dt{unknown} general nonlinear law $\Q(c_{1},v)$ that \dt{is autonomous in $c_{1}$ and $v$, which} will be reconstructed from the data available at a specific \dt{later} time \dt{in the tumour evolution}}. 
\end{enumerate}
%
%

\re{To complete the description of the model, we mention that the initial conditions for the cancer-ECM dynamics} described by equations \eqref{eq:Direct_Mutation_DepOn_C} are: 
 \bequ
 \label{eq:DP_IC}
 c_{1}(x,0):=c_{1,0}(x), \quad  c_{2}(x,0):=c_{2,0}(x) \quad \textrm{and} \quad v(x,0):=v_{0}(x), \;\; \textrm{for} \;\; x\in \Omega.
 \eequ
Here, $c_{1,0}(\cdot)$, $c_{2,0}(\cdot)$ and $v_{0}(\cdot)$ \re{give the initial distributions of the} primary cell subpopulation, mutated cell subpopulation and ECM, respectively. Furthermore, we assume that the cells \re{and the ECM components} do not leave the tissue region $\Omega$, \re{and therefore} we consider zero Neumann boundary conditions: 
\bequ\label{eq:DP_BC}
\frac{\partial c_{1}}{\partial n}\bigg |_{\partial \Omega}=0, \quad \frac{\partial c_{2}}{\partial n}\bigg |_{\partial \Omega}=0 \quad  \textrm{and} \quad \frac{\partial v}{\partial n}\bigg |_{\partial \Omega}=0, 
\eequ
where $n(\xi)$ is the usual normal direction at any given tissue boundary point $\xi\in\partial \Omega $.

\re{Throughout} the following sections, \re{we refer to} the tumour dynamics \re{model} \eqref{eq:Direct_Mutation_DepOn_C} together with the initial and boundary conditions \eqref{eq:DP_IC} and \eqref{eq:DP_BC} as the \emph{``forward model"}. 
%
%
%
%
%
%
%
%
\section{Inverse Problem Formulation for the Unknown Cancer Cell Mutation}
\label{Inverse_Problem_for_Unknown_Cancer_Cells_Mutation}
\re{We start with} the \emph{forward model} defined by the tumour dynamics \eqref{eq:Direct_Mutation_DepOn_C} in the presence of the initial and boundary conditions \eqref{eq:DP_IC} and \eqref{eq:DP_BC}. \re{Our goal is} to reconstruct the unknown cancer cells mutation law $\Q(\cdot)$ from additional information \re{given} by measurements of the cancer cells and ECM densities taken at \re{some} later time $t_{f}:=T>0$ in the tumour evolution. These measurements are therefore given in the form of \re{the following} functions on $\Omega$, \re{which} are considered to be known:
\begin{subequations}
\label{eq:Additional_Information_Two_Pop}
\begin{align}
  c^{*}_{1}(\cdot):\Omega\to \R &\quad \textrm{for the cancer subpopulation $c_{1}$,}
  \label{eq:Additional_Information_One_Pop_C1_Nonlocal}\\ 
  c^{*}_{2}(\cdot):\Omega\to \R &\quad \textrm{for the cancer subpopulation $c_{2}$,}
    \label{eq:Additional_Information_One_Pop_C2_Nonlocal}\\ 
 v^{*}(\cdot):\Omega\to \R &\quad \textrm{for the ECM density.} 
 \label{eq:Additional_Information_One_Pop_V_Nonlocal}
\end{align}
\end{subequations} 
In the following, we will explore the reconstruction of the unknown cancer mutation law $\Q(\cdot)$ when the known measurements $c^{*}_{1}(x)$, $c^{*}_{2}(x)$ and $v^{*}(x)$ will be given both as exact (accurate) data and as noisy data, $\forall x \in \Omega$. 
%
%
%
%
\subsection{Inverse Problem Setup: Forward Solver \dt{for the Retrieval of} \re{Mutation Laws} \dt{in Cases (i) and (ii) }}
\label{Inverse_Problem_Forward_Operator_and_Computational_Approach}
\dt{In this section we outline in a unitary manner the forward solver involved in the retrieval of the mutation laws corresponding to cases (i) and (ii) that require the retrieval of $\widetilde{\Q}_{1}(c_{1})$ and $\widetilde{\Q}_{2}(v)$, respectively. To that end, for $r=1,2$, denoting by $e_{r}$ either the primary tumour or the ECM, \emph{i.e.},  $e_{r}\in \{c_{1}, v\}$, enables us to proceed with addressing simultaneously both cases by simply referring to the retrieval of the term compactly denoted as $\widetilde{\Q}_{r}(e_{r})$ that is specified by
\bequ\label{unitary_notation_07oct21}
e_{r}:=
\left\{
\begin{array}{ll}
c_{1}, & if \quad r=1,\\ 
v, & if \quad r=2,
\end{array}
\right.
\quad
\textrm{and subsequently}\quad
\widetilde{\Q}_{r}(e_{r}):=
\left\{
\begin{array}{ll}
\widetilde{\Q}_{1}(c_{1}), & if \quad r=1,\\ 
\widetilde{\Q}_{2}(v), & if \quad r=2.
\end{array}
\right.
\eequ}
\dt{We start by} \re{considering} a uniform discretisation $\G_{_{\Omega}}:=\{(x_{i}, y_{j})\}_{i,j=1\dots N}$ of step size $\Delta x=\Delta y>0$ for a square maximal tissue region $\Omega\subset \R^{2}$ where the tumour exercises its dynamics. At any given time $t\in[0, t_{f}]$ the discretisations of cancer \dt{cells} densities $c_{1}\left(\cdot, t\right)$ \dt{and} $c_{2}\left(\cdot, t\right)$ \dt{as well as the density of ECM} $v\left(\cdot, t\right)$ are therefore given by the $N\times N$ matrices $\tilde c_{1}\left(t\right) := \{\tilde c_{1,i,j}\left(t\right)\}_{i,j=1\dots N}$, $\tilde c_{2}\left(t\right) := \{\tilde c_{2,i,j}\left(t\right)\}_{i,j=1\dots N}$  and $\tilde v\left(t\right) := \{\tilde v_{i,j}\left(t\right)\}_{i,j=1\dots N}$, with $\tilde c_{1,i,j}(t) := c_{1}((x_{i}, y_{j}), t)$, $\tilde c_{2,i,j}(t) := c_{2}((x_{i}, y_{j}), t)$ and $\tilde v_{i,j}(t) := v((x_{i}, y_{j}), t)$, $\forall \,\, i,j=1\dots N$. \dt{Correspondingly, in the following $\tilde{e}_{r}(t)$ will denote either $\tilde{c}_{1}(t)$ or $\tilde{v}(t)$, \emph{i.e.}, $\tilde{e}_{r}(t)\in\{\tilde{c}_{1}(t), \tilde{v}(t)\}$, as required by \eqref{unitary_notation_07oct21}}.

\re{Throughout this study we assume that we have \emph{a priori} knowledge that the cumulated ECM and cancer densities do not exceed the tissue carrying capacity $K_{c}$. Under this assumption, the unknown mutation law can} be written in terms of an unknown (for the moment) function \dt{$m^{c_{1}^{*}, c_{2}^{*}, v^{*}}:[0,K_{c}] \to \left[0, \infty\right)$}. Moreover, this unknown function \dt{$m^{c_{1}^{*},c_{2}^{*},v^{*}}$} will be appropriately identified within a suitable family of functions \dt{ $\Mu^{1}$} such that the corresponding solution for the tumour model \eqref{eq:Direct_Mutation_DepOn_C} will match the measurements given in \eqref{eq:Additional_Information_Two_Pop}. Thus, denoting by $\widetilde{\Q}_{\dt{r}}^{c_{1}^{*}, c_{2}^{*}, v^{*}}(\cdot) $ the unknown mutation \dt{term} for which the corresponding solution of model \eqref{eq:Direct_Mutation_DepOn_C} matches measurement \eqref{eq:Additional_Information_Two_Pop}, at each $(x_{i}, y_{j})$ we can write this as 
\bequd
\dt{\widetilde{\Q}_{\dt{r}}^{c_{1}^{*}, c_{2}^{*}, v^{*}}(\tilde{e}^{m^{c_{1}^{*},c_{2}^{*},v^{*}}}_{r,i,j}(t)):= \M^{1}_{i,j}(\tilde{e}_{r}^{m^{c_{1}^{*},c_{2}^{*},v^{*}}}(t), m^{c_{1}^{*},c_{2}^{*},v^{*}}),}
\eequd
where $\dt{\M^{1}\left(\cdot, \cdot\right) := \{\M^{1}_{i,j}}(\cdot, \cdot)\}_{i,j=1\dots N}$, \dt{with} $\M^{1}(\cdot, \cdot) : \R^{N \times N} \times \dt{\Mu^{1}}\to \R^{N\times N} $ \dt{representing} a \emph{``trial mutation operator"} that will be specified below alongside the family of functions \dt{$\Mu^{1}$. Furthermore, $\tilde{e}_{r}^{m^{c_{1}^{*},c_{2}^{*},v^{*}}}(t):=\{\tilde{e}^{m^{c_{1}^{*},c_{2}^{*},v^{*}}}_{r,i,j}(t)\}_{i,j=1..N}$ represents the solution for the density of either the primary cell population (if $r=1$) or the ECM (if $r=2$) that is obtained for model \eqref{eq:Direct_Mutation_DepOn_C} when, instead of the unknown term $\widetilde{\Q}_{r}(e_{r})$, in the mutation law we use the trial mutation term $\M^{1}(\cdot, m^{c_{1}^{*},c_{2}^{*},v^{*}})$.} 

\re{Next, we consider} an uniform discretisation for the domain $\left[0, K_{c}\right]$ that is given by an equally spaced grid $\dt{\G^{1}_{_{M}}} := \left\{\eta_{l}\right\}_{l=1\dots M}$ of step size $\Delta \eta>0$. \re{On this discretised domain, the} unknown function \dt{$m^{c_{1}^{*},c_{2}^{*},v^{*}}$ is} identified through a suitable approximation within the following $M-$dimensional space of functions associated with \dt{$\G^{1}_{_{M}}$}, namely 
\bequ
\begin{split}
\dt{\Mu^{1}} := &\bigg\{m : \left[0, K_{c}\right] \to \R\,\bigg |\,\, m\restrict{E_{l}} = \sum\limits_{p=0,1} m(\eta_{l+p}) \phi_{l+p},\quad \forall E_{l}\in \G_{_{M}}^{1,tiles}\bigg\} \\
&\textrm{with the family of intervals} \;\; \dt{\G_{_{M}}^{1,\,tiles}}:=\left\{ E_{l} := \left[\eta_{l}, \eta_{l+1}\right]\,|\, l=1\dots M-1\right\},  \\
& \text{and} \;\; \forall \,\, E_{l}\in \dt{\G_{_{M}}^{1,\,tiles}}, \,\,\{\phi_{l+p}\}_{p=0,1}\; \textrm{\re{describe} the usual linear shape functions on $E_{l}$.} 
\end{split}
\eequ
Thus, for any candidate function $m\in \dt{\Mu^{1}}$, the corresponding  \emph{trial mutation operator} \dt{$\M^{1}$} has each of its components $\M_{i,j}$, $\forall \,i,j=1\dots N$, given by
\bequ
\label{attempted_mutation}
\begin{split}
&\dt{\qquad\qquad \M^{1}_{i,j}(\tilde{e}_{r}^{m}(t), m) := m\restrict{E_{l}}(\tilde{e}^{m}_{r,i,j}(t)),} \\
& \textrm{with index $l$ being independent of its choice within the associated set of indices $\Lambda_{i,j}$, namely:}\\
&\Lambda_{i,j} := \{l\in \{1, \dots, M\!-\!1\}\,\,|\,\, \exists E_{l}\in \dt{\G_{_{M}}^{1,\,tiles}} \textrm{ such that } \tilde{e}^{m}_{r,i,j}(t)\in E_{l} \}. 
\end{split}
\eequ
{\color{blue}Here, \dt{as per \eqref{unitary_notation_07oct21}}}, \dt{$\tilde{e}_{r}^{m}(t) := \{\tilde{e}^{m}_{ri,j}(t) \}_{i,j=1\dots N}$ represents either the solutions for the density of primary cancer cell population, $\{\tilde{c}^{m}_{1,i,j}(t) \}_{i,j=1\dots N}$,  or for the density of ECM, $\{\tilde{v}^{m}_{i,j}(t) \}_{i,j=1\dots N}$, which alongside the density of mutated cell population, $\{\tilde{c}^{m}_{1,i,j}(t) \}_{i,j=1\dots N}$, are obtained with model \eqref{eq:Direct_Mutation_DepOn_C} when this uses within the mutation law the trial mutation term $\M^{1}(\cdot,m):=\{\M^{1}_{i,j}\left(\cdot, m\right)\}_{i,j=1\dots N}$ given in \eqref{attempted_mutation} instead of the unknown term $\widetilde{\Q}_{r}(\cdot)$. Finally, the \emph{trial mutation form} for the full mutation law given in \eqref{eq:Mutation_DepOn_C} and \eqref{eq:Mutation_DepOn_V} for cases $(i)$ and $(ii)$, respectively, is denoted by $\overline{\M}_{r}(\tilde{c}^{m}_{1}(t),\tilde{v}^{m}(t),m):=\{\overline{\M}_{r,i,j}(\tilde{c}^{m}_{1}(t),\tilde{v}^{m}(t),m)\}_{i.j=1\dots N}$ and is given by:
\bequ
\overline{\M}_{r,i,j}(\tilde{c}^{m}_{1}(t),\tilde{v}^{m}(t),m):=
\left\{
\begin{array}{ll}
\M^{1}_{i,j}(\tilde{c}_{1}^{m}(t), m), & if \quad r=1,\\[0.2cm] 
\delta_{0}\tilde{c}_{1,i,j}^{m}(t)\M^{1}_{i,j}(\tilde{v}^{m}(t), m), & if \quad r=2.
\end{array}
\right.
\eequ
Therefore, in space-discretised form, model \eqref{attempted_mutation} that uses $\M^{1}\left(\cdot, m\right)$ can therefore be written as}
\bequ
\label{eq:spatial_discretization_1}
\frac{\partial}{\partial{t}}
\left[
\begin{array}{l}
\tilde c^{m}_{1} \\
\tilde c^{m}_{2} \\
\tilde v^{m}
\end{array}
\right]
= 
\left[
\begin{array}{l}
\H^{1}(\tilde c_{1}^{m}, c_{2}^{m}, \tilde{v}^{m}, m)\\
\H^{2}(\tilde c_{1}^{m}, c_{2}^{m}, \tilde{v}^{m}, m)\\
\H^{3}(\tilde c_{1}^{m}, c_{2}^{m}, \tilde{v}^{m})\\
\end{array}
\right], 
\eequ
Here, $\H^{1}\left(\cdot, \cdot, \cdot, \cdot\right) = \{\H^{1}_{i,j}\left(\cdot, \cdot, \cdot, \cdot\right)\}_{i,j=1\dots N}$ represents the spatial discretisation corresponding to equation \eqref{eq:DP_C1_Simp_Mut_Chap4}. Each of its components $\H^{1}_{i,j}\left(\cdot, \cdot, \cdot, \cdot\right)$, $\forall\,\, i,j=1\dots N$, are given by 
\bequ
\label{opH1}
\begin{split}
\H^{1}_{i,j} & \left(\tilde{c}_{1}^{m}\left(t\right), \tilde{c}_{2}^{m}\left(t\right), \tilde{v}^{m}\left(t\right), m\right) := \\
& \frac{D_{1}}{\left(\Delta{x}\right)^{2}}\left(\tilde{c}^{m}_{1,i-1,j}\left(t\right) + \tilde{c}^{m}_{1,i+1,j}\left(t\right) + \tilde{c}^{m}_{1,i,j-1}\left(t\right) + \tilde{c}^{m}_{1,i,j+1}\left(t\right) - 4 \tilde{c}^{m}_{1,i,j}\left(t\right) \right) \\
& - \frac{\eta_{1}}{2\left(\Delta{x}\right)^{2}} \left(\left(\tilde{c}^{m}_{1,i,j}\left(t\right) + \tilde{c}^{m}_{1,i+1,j}\left(t\right)\right)\left(\tilde{v}^{m}_{i+1,j}\left(t\right) - \tilde{v}^{m}_{i,j}\left(t\right)\right) - \left(\tilde{c}^{m}_{1,i,j}\left(t\right) + \tilde{c}^{m}_{1,i-1,j}\left(t\right)\right)\left(\tilde{v}^{m}_{i,j}\left(t\right) - \tilde{v}^{m}_{i-1,j}\left(t\right)\right)\right. \\
& + \left.  \left(\tilde{c}^{m}_{1,i,j}\left(t\right) + \tilde{c}^{m}_{1,i,j+1}\left(t\right)\right)\left(\tilde{v}^{m}_{i,j+1}\left(t\right) - \tilde{v}^{m}_{i,j}\left(t\right)\right) - \left(\tilde{c}^{m}_{1,i,j}\left(t\right) + \tilde{c}^{m}_{1,i,j-1}\left(t\right)\right)\left(\tilde{v}^{m}_{i,j}\left(t\right) - \tilde{v}^{m}_{i,j-1}\left(t\right)\right)\right) \\
& + \mu_{c} \tilde c^{m}_{1,i,j}\left(t\right) \left(K_{c} - \tilde{c}^{m}_{1,i,j}\left(t\right) - \tilde{c}^{m}_{2,i,j}\left(t\right) - \tilde{v}^{m}_{i,j}\left(t\right)\right) - \dt{\omega(t)\overline{\M}_{r,i,j}(\tilde{c}^{m}_{1}(t),\tilde{v}^{m}(t),m)}. \\
\end{split}
\eequ
%
Similarly, $\H^{2}\left(\cdot, \cdot, \cdot, \cdot\right) = \{\H^{2}_{i,j}\left(\cdot, \cdot, \cdot, \cdot\right)\}_{i,j=1\dots N}$ represents the spatial discretisation corresponding to equation \eqref{eq:DP_C2_Simp_Mut_Chap4}, and each of its components $\H^{2}_{i,j}\left(\cdot, \cdot, \cdot, \cdot\right)$, $\forall\,\, i,j=1\dots N$, are given by 
\bequ
\label{opH2}
\begin{split}
\H^{2}_{i,j}&(\tilde{c}_{2}^{m}(t), \tilde{c}_{2}^{m}(t), \tilde{v}^{m}(t), m) :=\\
&\frac{D_{2}}{\left(\Delta{x}\right)^{2}}\left(\tilde{c}^{m}_{2,i-1,j}(t) + \tilde{c}^{m}_{2,i+1,j}(t) + \tilde{c}^{m}_{2,i,j-1}(t) + \tilde{c}^{m}_{2,i,j+1}(t) - 4 \tilde{c}^{m}_{2,i,j}(t) \right) \\
& - \frac{\eta_{2}}{2\left(\Delta{x}\right)^{2}} \left(\left(\tilde{c}^{m}_{2,i,j}(t) + \tilde{c}^{m}_{2,i+1,j}(t)\right)\left(\tilde{v}^{m}_{i+1,j}(t) - \tilde{v}^{m}_{i,j}(t)\right) - \left(\tilde{c}^{m}_{2,i,j}(t) + \tilde{c}^{m}_{2,i-1,j}(t)\right)\left(\tilde{v}^{m}_{i,j}(t) - \tilde{v}^{m}_{i-1,j}(t)\right)\right. \\
& + \left.  \left(\tilde{c}^{m}_{2,i,j}(t) + \tilde{c}^{m}_{2,i,j+1}(t)\right)\left(\tilde{v}^{m}_{i,j+1}(t) - \tilde{v}^{m}_{i,j}(t)\right) - \left(\tilde{c}^{m}_{2,i,j}(t) + \tilde{c}^{m}_{2,i,j-1}(t)\right)\left(\tilde{v}^{m}_{i,j}(t) - \tilde{v}^{m}_{i,j-1}(t)\right)\right) \\
& + \mu_{c} \tilde c^{m}_{2,i,j}(t) (K_{c} - \tilde{c}^{m}_{1,i,j}(t) - \tilde{c}^{m}_{2,i,j}(t) - \tilde{v}^{m}_{i,j}(t)) + \dt{\omega(t)\overline{\M}_{r,i,j}(\tilde{c}^{m}_{1}(t),\tilde{v}^{m}(t),m)}. \\
\end{split}
\eequ
\re{Finally,}  $\H^{3}(\cdot, \cdot, \cdot)=\{\H^{3}_{i,j}(\cdot,\cdot, \cdot)\}_{i,j=1\dots N}$ represents the discretisation of the ECM equation \eqref{eq:DP_V_Simp_Mut_Chap4}, and each of its components $\H^{3}_{i,j}(\cdot,\cdot)$, $\forall\,\, i,j=1\dots N$, are given by
\bequ
\label{opH3}
\H^{3}_{i,j}(\tilde{c}_{1}^{m}(t), \tilde{c}_{2}^{m}(t), \tilde{v}^{m}(t)) := - \rho(\tilde{c}_{1,i,j}^{m}(t) + \tilde{c}_{2,i,j}^{m}(t))\tilde{v}_{i,j}^{m}(t)+ \mu_{v}(K_{c} - \tilde{c}_{1,i,j}^{m}(t) - \tilde{c}_{2,i,j}^{m}(t) -\tilde{v}_{i,j}^{m}(t))^{+}.
\eequ
\re{Consider now} a time discretisation $\{t_{n}\}_{n=0\dots L}$ with time step $\Delta t:=T/L$. For each $n\in\{\dt{0},\dots,L\}$, a simple Euler time-marching scheme can be \re{written} for system \eqref{eq:spatial_discretization_1} via the associated operator 
\bequ\label{MOL_1}
\begin{split}
&K_{m}: \R^{N\times N}\times \R^{N\times N} \times \R^{N\times N} \to \R^{N\times N}\times \R^{N\times N} \times \R^{N\times N}\\
\textrm{given by}&\\
&K_{m}
\left(
\begin{bmatrix}
\tilde c^{m,n}_{1}\\ 
\tilde c^{m,n}_{2}\\ 
\tilde v^{m,n}
\end{bmatrix}
\right)
:= 
\begin{bmatrix}
\tilde c^{m,n}_{1}\\ 
\tilde c^{m,n}_{2}\\ 
\tilde v^{m,n}
\end{bmatrix}
+ \Delta{t} 
\begin{bmatrix}
\H^{1}(\tilde c^{m, n}_{1}, \tilde c^{m, n}_{2}, \tilde{v}^{m,n}, m)\\
\H^{2}(\tilde c^{m, n}_{1}, \tilde c^{m, n}_{2}, \tilde{v}^{m,n}, m)\\
\H^{3}(\tilde c^{m, n}_{1}, \tilde c^{m, n}_{2}, \tilde{v}^{m,n})
\end{bmatrix},
\end{split}
\eequ
where $\tilde c^{m,n}_{1}:=\tilde c^{m}_{1}(t_{n})$, $\tilde c^{m,n}_{2}:=\tilde c^{m}_{2}(t_{n})$ and $\tilde{v}^{m,n}:=\tilde{v}^{m}(t_{n})$. The right-hand-side operators are given by  $\H^{1}(\tilde{c}^{m,  n}_{1}, \tilde{c}^{m,  n}_{2}, \tilde{v}^{m,n},m):=\H^{1}(\tilde c^{m}_{1}(t_{n}), \tilde c^{m}_{2}(t_{n}), \tilde{v}^{m}(t_{n}),m)$, $\H^{2}(\tilde{c}^{m,  n}_{1}, \tilde{c}^{m,  n}_{2}, \tilde{v}^{m,n},m):=\H^{2}(\tilde c^{m}_{1}(t_{n}), \tilde c^{m}_{2}(t_{n}), \tilde{v}^{m}(t_{n}),m)$ and $\H^{3}(\tilde{c}^{m,  n}_{1}, \tilde{c}^{m,  n}_{2}, \tilde{v}^{m,n}):=\H^{3}(\tilde c^{m}_{1}(t_{n}), \tilde c^{m}_{2}(t_{n}), \tilde{v}^{m}(t_{n}))$. This \re{allows} us to formulate \re{the} \emph{``forward operator" $K$} between the family of functions \dt{$\Mu^{1}$} (where we search for the appropriate cancer cells mutation function \dt{$m^{c_{1}^{*}, c_{2}^{*}, v^{*}}$}) and the space where the discretised measurements \eqref{eq:Additional_Information_Two_Pop} are recorded. Hence, the \emph{``forward operator" $K$} is defined as
\bequ\label{MOL_2}
\begin{split}
&K:\S\to \R^{N\times N}\times \R^{N\times N} \times \R^{N\times N}\\
\textrm{given by}&\\
&K(m):=\underbrace{K_{m} \circ K_{m} \circ  \cdots \cdots  \circ K_{m} }_{\,\,\,\,\,\,L\,\,\, \text{times}}\left(
\begin{bmatrix}
\tilde{c}_{1,0}\\ 
\tilde{c}_{2,0}\\ 
\tilde{v}_{0}
\end{bmatrix}
\right)
\end{split}
\eequ
where $\tilde{c}_{1,0}:=\{c_{1,0}(x_{i},y_{j})\}_{i,j=1,\dots, N}$, $\tilde{c}_{2,0}:=\{c_{2,0}(x_{i},y_{j})\}_{i,j=1,\dots, N}$  and $\tilde{v}_{0}:=\{v_{0}(x_{i},y_{j})\}_{i,j=1,\dots, N}$ are the discretised initial conditions \eqref{eq:DP_IC} for the governing tumour \emph{forward model} \dt{\eqref{eq:Direct_Mutation_DepOn_C}}. Hence, for each \dt{$m \in \Mu^{1} $}, the forward operator $K$ gives the spatio-temporal progression of the initial condition $[\tilde{c}_{1,0}, \tilde{c}_{2,0}, \tilde{v}_{0}]^{T}$ under the invasion model \eqref{eq:Direct_Mutation_DepOn_C}, which is obtained when the cell mutation law at each instance of time $t>0$ \dt{involves} the \dt{corresponding} trial mutation operator \dt{$\M(\cdot,m)$ instead of the unknown mutation terms $\widetilde{\Q}_{r}(\cdot)$}. 
\subsection{The Inverse Problem Regularisation Approach \re{for Mutation Laws in Cases (i) and (ii) }}
\label{Numerical_Solution_of_Inverse_Problem}
From \eqref{MOL_1} and \eqref{MOL_2} we have that our forward operator $K$ is given as a finite composition of affine functions of the form 
\bequ
\dt{\Mu^{1}}\ni m\longmapsto K_{m}\in \boldsymbol{\ell}^{2}(\boldsymbol{\ell}^{2}(\E\times \E\times \E); \boldsymbol{\ell}^{2}(\E\times\E\times \E)).
\eequ
\re{Here,} $\boldsymbol{\ell}^{2}(\boldsymbol{\ell}^{2}(\E\times \E\times \E); \boldsymbol{\ell}^{2}(\E\times \E\times \E))$ \re{is} the usual finite-dimensional Bochner space of square integrable vector-value functions \cite{yosida1980} with respect to the counting measure (see \cite{Schilling2005}, p. 27) that are defined on $\boldsymbol{\ell}^{2}(\E\times \E\times \E)$ and take values in $\boldsymbol{\ell}^{2}(\E\times \E\times \E)$, \re{and} $\E:=\{E_{i,j}\}_{i,j=1\dots N}$ represents the standard basis of elementary matrices associated with the grid $\G_{_{\Omega}}$.  As a direct consequence, we immediately obtain that this \re{operator}  is both continuous and compact, from where we obtain that $K$ is also closed and sequentially bounded \cite{yosida1980}. Therefore, $K$ satisfies the hypotheses assumed in \cite{Engl_1989} that ensure convergence for the nonlinear Tikhonov regularisation strategy given by the functionals $\{J_{\alpha}\}_{\alpha>0}$,
\bequ
\label{eq:Objective_Function_1}
\begin{split}
&J_{\alpha}:\dt{\Mu^{1}}\to\R, \quad \forall \alpha>0,\\
\textrm{defined by}&\\
&J_{\alpha}\left(m\right) := \left\|K(m)- \begin{bmatrix}
\tilde{c}^{*}_{1}\\ 
\tilde{c}^{*}_{2}\\ 
\tilde{v}^{*}
\end{bmatrix}\right\|^{2}_{2} + \alpha \|{m}\|^{2}_{2},\quad \forall m\in \dt{\Mu^{1}}.
\end{split}
\eequ
\re{The minimisation of these functionals} enable us to identify \dt{$m^{c_{1}^{*}, c_{2}^{*}, v^{*}}$} as the limit $\alpha\to 0$ of the points of minimum $m^{\alpha}$ of $J_{\alpha}$ \re{(these points correspond to the smallest discrepancy between the data measurements and the solution of our system that uses $m^{\alpha}$ as a mutation law)}. The two norms involved in \eqref{eq:Objective_Function_1} represent the usual Euclidean norms on the corresponding finite dimensional spaces. Indeed, while the first is the standard  Euclidean norm on $\R^{N\times N}\times \R^{N\times N} \times \R^{N\times N}$, the second is also the Euclidean norm induced on the $M-$dimensional space of functions $\Mu^{1}$ via the standard isomorphism that we establish between \dt{$\Mu^{1}$} and $\R^{M}$ by which each $m\in \dt{\Mu^{1}}$ is uniquely represented through its nodal values $\{m(\eta_{l})\}_{l=1\dots M}$ with respect to the \dt{linear} basis functions $\{\bar \phi_{l}\}_{l=1\dots M}$ associated to $\G_{M}$ \cite{Hughes1987}: 
\bequ
\begin{split}
\textrm{since}\,\,m=\!\!\sum\limits_{l=1\dots M} \!\!\!m(\eta_{l})\bar \phi_{l},
\,\, \textrm{ we therefore make the identification: }\,\,
m\equiv \{m(\eta_{l})\}_{l=1\dots M}.\\
\end{split}
\eequ
Finally, in \eqref{eq:Objective_Function_1}, $\tilde{c}^{*}_{1}$, $\tilde{c}^{*}_{2}$ and $\tilde{v}^{*}$ represent the discretised measurements of the densities of cancer cells and ECM given in equations \eqref{eq:Additional_Information_One_Pop_C1_Nonlocal}-\eqref{eq:Additional_Information_One_Pop_V_Nonlocal}, i.e., $\tilde{c}^{*}_{1}:=\{c^{*}_{1}(x_{i},y_{j})\}_{i,j=1,\dots, N}$, $\tilde{c}^{*}_{2}:=\{c^{*}_{2}(x_{i},y_{j})\}_{i,j=1,\dots, N}$  and $\tilde{v}^{*}:=\{v^{*}(x_{i},y_{j})\}_{i,j=1,\dots, N}$.
We \re{assume that these data} measurements are either exact or are corrupted by a certain \re{noise} level $\delta\geq 0$. Thus, maintaining for simplicity the measurements notation unchanged, these \re{measurements} are given by
\begin{subequations}
\label{comp_meas_1}
\begin{align}
\tilde{c}^{*}_{1}(x) &= \tilde{c}^{*}_{1exact}(x) + \delta \gamma_{c_{1}}(x),\label{comp_meas_c1_1}\\
\tilde{c}^{*}_{2}(x) &= \tilde{c}^{*}_{2exact}(x) + \delta \gamma_{c_{2}}(x),\label{comp_meas_c2_1}\\
\tilde{v}^{*}(x) &= \tilde{v}^{*}_{exact}(x) + \delta \gamma_{v}(x), \label{comp_meas_v_1}
\end{align}
\end{subequations}
where, $\forall\,\, x\in \Omega$, we have that $\tilde{c}^{*}_{1exact}(x)$, $\tilde{c}^{*}_{2exact}(x)$ and $\tilde{v}^{*}_{exact}(x)$ \re{describe} the exact data, and $\gamma_{c_{1}}(x)$, $\gamma_{c_{2}}(x)$ and  $\gamma_{v}(x)$ are signal-independent noise generated from a Gaussian normal distribution with mean zero and standard deviations $\sigma_{c_{1}}$, $\sigma_{c_{2}}$ and $\sigma_{v}$, respectively, given by
\begin{equation}
\label{noise_in_the_measurements}
\begin{cases}
\sigma_{c_{1}} := \frac{1}{\lambda \left(\Omega\right)} \int\limits_{\Omega}\tilde{c}^{*}_{1exact}(x) \,dx,\\
\sigma_{c_{2}} := \frac{1}{\lambda \left(\Omega\right)} \int\limits_{\Omega}\tilde{c}^{*}_{2exact}(x) \,dx,\\
\sigma_{v} := \frac{1}{\lambda \left(\Omega\right)} \int\limits_{\Omega}\tilde{v}^{*}_{exact}(x) \,dx,
\end{cases}
\end{equation}
with $\lambda\left(\cdot\right)$ being the usual Lebesgue measure. In the numercial results below, we generate the random variables $\gamma_{c_{1}}(x)$, $\gamma_{c_{2}}(x)$ and $\gamma_{v}(x)$ via MATLAB function \textit{normrnd} by taking $\{\gamma_{c_{1}}(x_{i},y_{j})\}_{i,j=1\dots N}:= \textit{normrnd}\left(0, \sigma_{c_{1}}, N \times N\right)$, $\{\gamma_{c_{2}}(x_{i},y_{j})\}_{i,j=1\dots N}:= \textit{normrnd}\left(0, \sigma_{c_{2}}, N \times N\right)$ and $\{\gamma_{v}(x_{i},y_{j})\}_{i,j=1\dots N}:= \textit{normrnd}\left(0, \sigma_{v}, N \times N\right)$.
%
%
%
%
%
%
\subsection{Numerical Reconstruction of the Unknown Mutation Laws Terms in cases (i) and (ii)}
\label{Minimization_Algorithms}

We explore now the inversion approach that we formulated so far in the context of \re{the} forward model \eqref{eq:Direct_Mutation_DepOn_C} by proceeding with the reconstruction of \dt{the unknown terms involved in the mutation laws in cases (i) and (ii), namely $\widetilde{\Q}_{1}(c_{1})$ and $\widetilde{\Q}_{2}(v)$, respectively. }

\paragraph{\emph{Initial Conditions}} The initial conditions \eqref{eq:DP_IC} that we consider in the computations for the forward model \eqref{eq:Direct_Mutation_DepOn_C} are \re{as follows:}
 \begin{subequations}
 \label{subeq:Initial_Conditions_One_Pop_1}
\begin{align}
\label{eq:initial_condition_C1_1}
c_{1,0}(x) & := 0.5 \left( exp\left(-\frac{\|x - \left(2,2\right)\|^{2}_{2}}{0.03}\right) - exp\left( - 9.407\right) \right), \\
\label{eq:initial_condition_C2_1}
c_{2,0}(x) & := 0, \\
\label{eq:heterogenous_initial_condition_1}
v_{0}(x)  & := 0.5 + 0.3 \cdot sin\left(4\pi \cdot \|x\|_{2}\right),    \quad \quad \forall \,x \in{\Omega}.
\end{align}
\end{subequations}
Here, \re{we assume that} $c_{2,0}(x)=0$ because \re{this second cell population} will arise after a period of time \re{following} mutations of the first \re{cell} population $c_{1}$. 

To identify the cancer cells mutation law, we consider both exact and noisy measurement data \eqref{comp_meas_1} as additional information for the forward model \eqref{eq:Direct_Mutation_DepOn_C} in the presence of initial conditions \eqref{subeq:Initial_Conditions_One_Pop_1} and boundary conditions \eqref{eq:DP_BC}. Specifically, \dt{in each of the two cases}, we \re{assume} that the exact data (namely $\tilde{c}^{*}_{1,exact}(x)$, $\tilde{c}^{*}_{2,exact}(x)$ and $\tilde{v}^{*}_{exact}(x)$) that appear in \eqref{comp_meas_1} are given by the solutions \dt{densities for primary tumour cells population, $\bar{c}_{1}(x,t)$, mutated cells population, $\bar{c}_{2}(x,t)$, and ECM, and $\bar v(x,t)$, evaluated} at the final time $t_{f}>0$, \emph{i.e.}, 
\bequ
 \label{exactMeasurementsOnePopulation}
\tilde{c}^{*}_{1,exact}(x):= \bar{c}_{1}(x,t_{f}),\quad \tilde{c}^{*}_{2,exact}(x):= \bar{c}_{2}(x,t_{f})\quad \textrm{and} \quad \tilde{v}^{*}_{exact}(x):=\bar{v}(x,t_{f}), \quad \forall x\in \Omega,
\eequ
\dt{which are obtained from the forward model \eqref{eq:Direct_Mutation_DepOn_C} as follows:
\begin{itemize}
\item for case $(i)$: we assume that the mutation law is of the form given in \eqref{eq:Mutation_DepOn_C}  but when parameter $\delta_{0}>0$ is considered known and has the value given in parameter Table \ref{paramSetTable}; \\[-0.1cm]
\item for case $(ii)$: we assume that the mutation law is of the form given by \eqref{eq:Mutation_DepOn_V} with the known term $\widetilde{\Q}_{2}(v)$ specified in \eqref{explicit_measurement_case_ii}.
\end{itemize}}
%
%

%
%
%
%
\begin{figure}[!ht]
\centering
\includegraphics[width=0.45\textwidth]{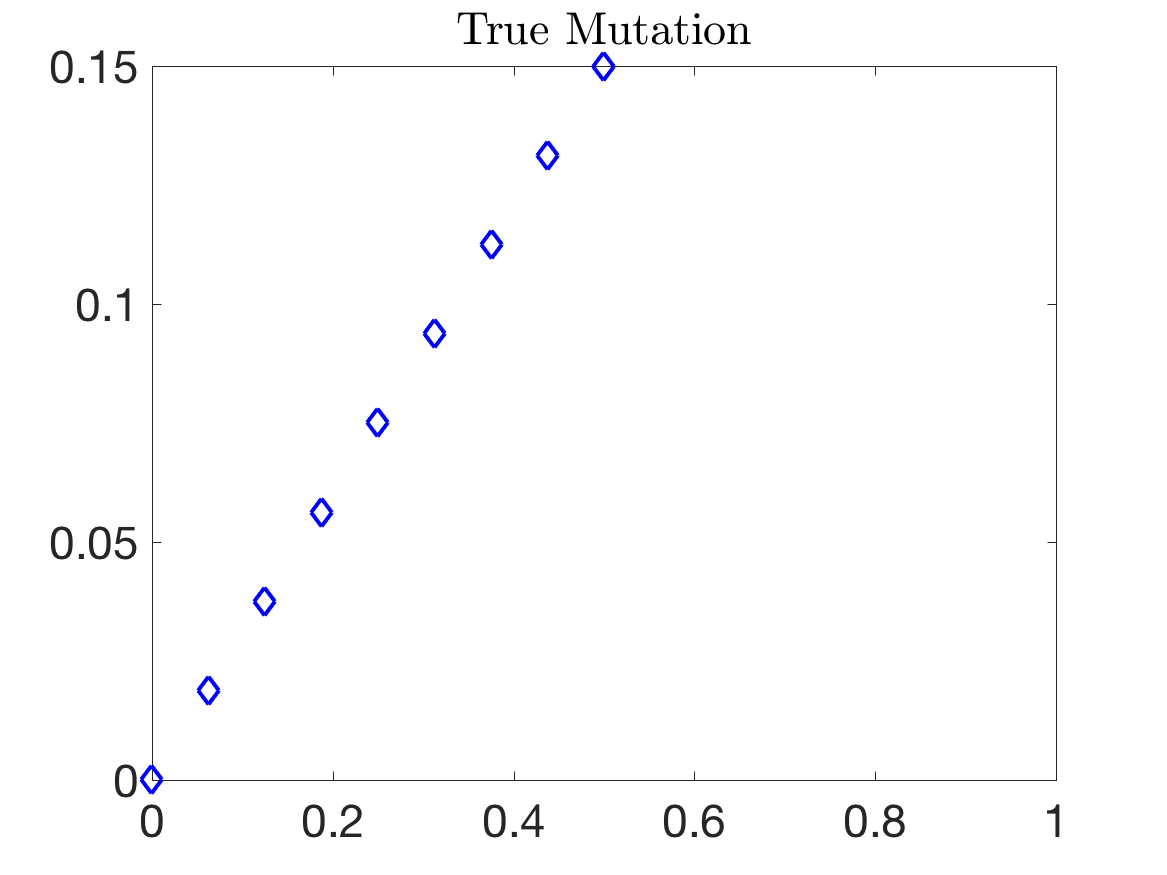}
\includegraphics[width=0.45\textwidth]{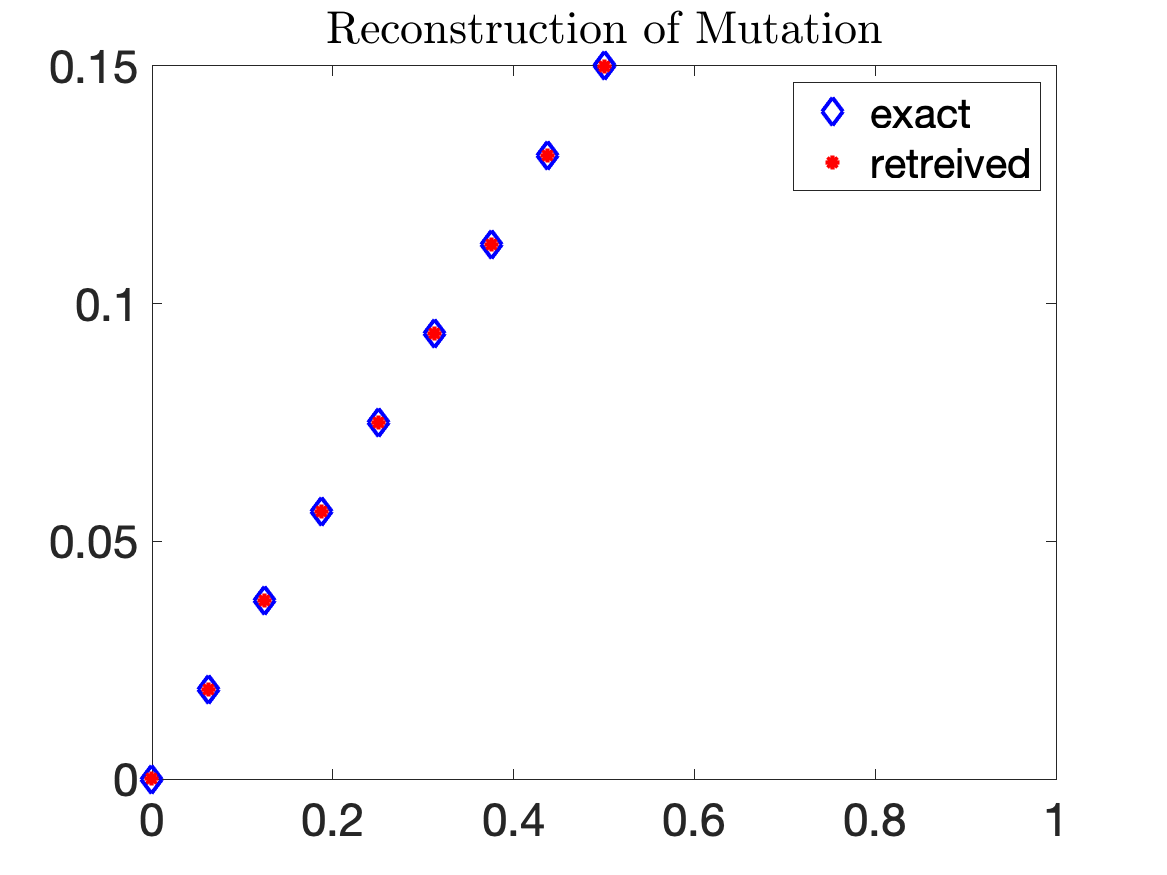}\\
\vspace{-0.1cm}
\includegraphics[width=0.45\textwidth]{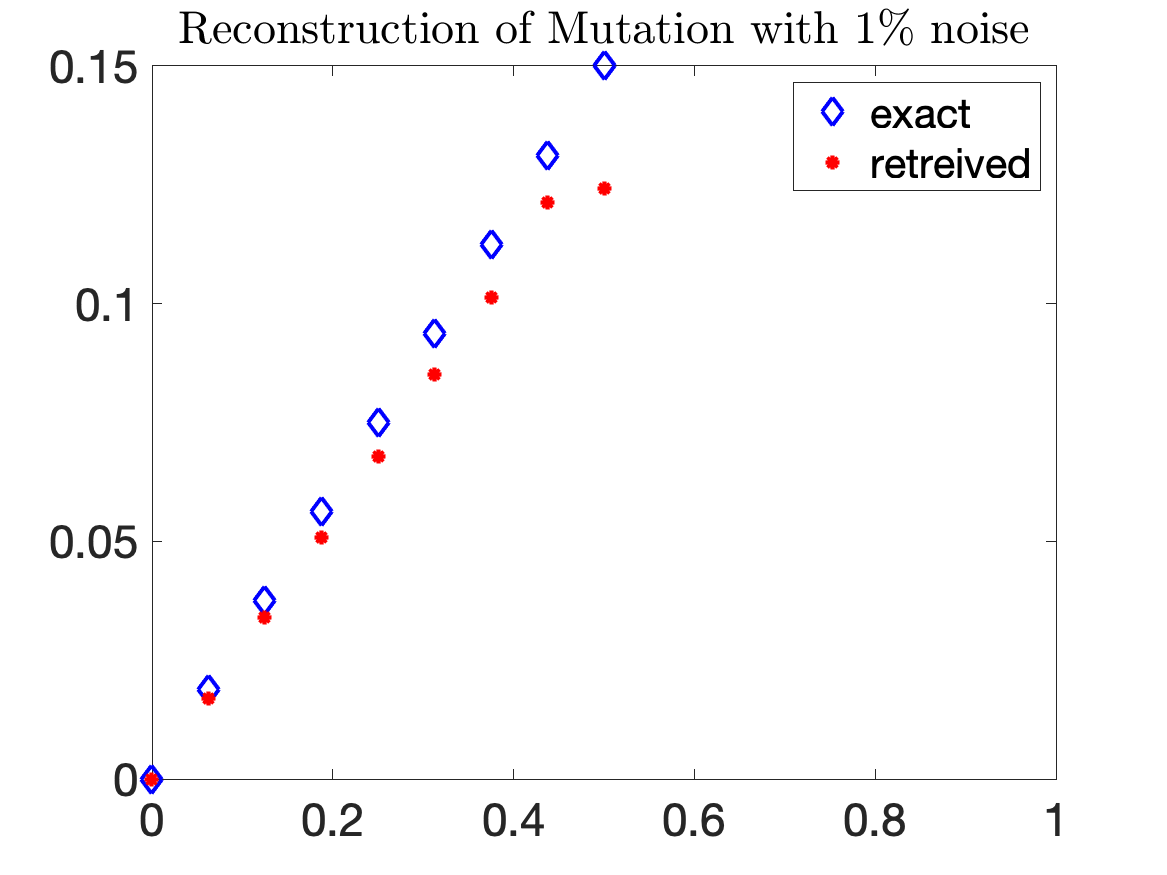}
\includegraphics[width=0.45\textwidth]{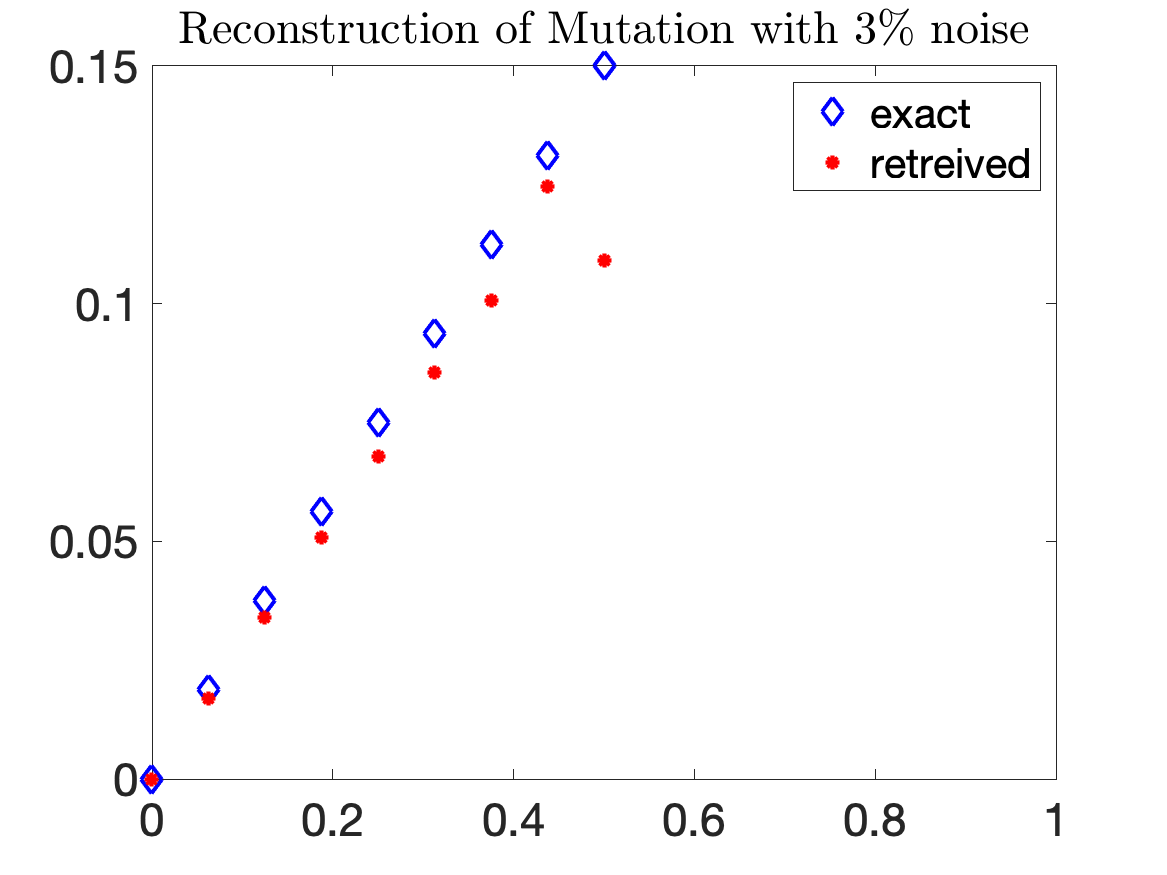}
\caption{\small{\textit{Reconstruction of \dt{the unknown} mutation law \re{\dt{term $\widetilde{\Q}_{1}(c_{1})$} \dt{that appears in mutation law given in \eqref{eq:Mutation_DepOn_C}} for} model \eqref{eq:Direct_Mutation_DepOn_C}. \re{First row: (left)} the true mutation law, and \re{(right)} the reconstructed mutation law restricted to $\A_{c}$ with no noise and $\alpha^{*}=10^{-12}$. \re{Second row: (left)} the reconstructed mutation law obtained for $1\%$ noisy data and $\alpha^{*}=10^{-5}$; \re{and (right) the reconstructed mutation law obtained} for $3\%$ noisy data and $\alpha^{*}=10^{-12}$. For all plots in this figure: 1) the first axis represents the the values for $c_{1}\in[\bar{c}_{1}^{min},\bar{c}_{1}^{max}]$; 2) second axis represents the values for mutation. \re{The numerical simulations are obtained using the parameters given in Table \ref{paramSetTable}.} }}}
\label{fig:True_Recon_Mut_Term_DepOn_C}
\end{figure}
For each regularisation parameter $\alpha>0$ considered here, the minimisation process for $J_{\alpha}$ is initiated with $m_{0} = \mathbb{I} \times10^{-3}$, (where $\mathbb{I}$ represents the $M$ vector of ones), and for the actual minimisation we employed here the nonlinear minimisation MATLAB function \textit{lsqnonlin}.
Finally, since there are no data to test the trial mutation operators beyond the \emph{maximal accessible region} $\A_{c}$ and $\A_{v}$ defined by the minimum and maximum values of the solution, i.e., 
\begin{equation}
\textrm{$\A_{c}:=[\bar{c}_{1}^{min},\bar{c}_{1}^{max}]$, with:} \quad
\bar{c}_{1}^{min} := \min\limits_{\begin{array}{c} \scriptstyle (x, t)\in \Omega \times [0,T]\\[-4pt] \end{array}} c_{1}(x, t), \quad \quad \quad \bar{c}_{1}^{max} := \max\limits_{\begin{array}{c} \scriptstyle (x, t)\in \Omega \times [0,T] \end{array}} c_{1}(x, t),
\label{EqAc}
\end{equation}
\begin{equation}
\textrm{$\A_{v}:=[\bar{v}^{min},\bar{v}^{max}]$, with:} \quad 
\bar{v}^{min} := \min\limits_{\begin{array}{c} \scriptstyle (x, t)\in \Omega \times [0,T]\\[-4pt] \end{array}} v(x, t), \quad \quad \quad \bar{v}^{max} := \max\limits_{\begin{array}{c} \scriptstyle (x, t) \in \Omega \times [0,T]\\[-4pt] \end{array}} v(x, t),\\
\label{EqAv}
\end{equation}
the reconstructions in this section will be attempted only for the restriction of the sought mutation laws to $\A_{c}$ and $\A_{v}$. \re{We need to emphasise that $\A_{c}$ is used to reconstruct \dt{$\widetilde{\Q}_{1}(c_{1})$}, while $\A_{v}$ is used to reconstruct \dt{$\widetilde{\Q}_{2}(v)$}. }

\begin{figure}[!ht]
	\centering
	\includegraphics[width=0.45\textwidth]{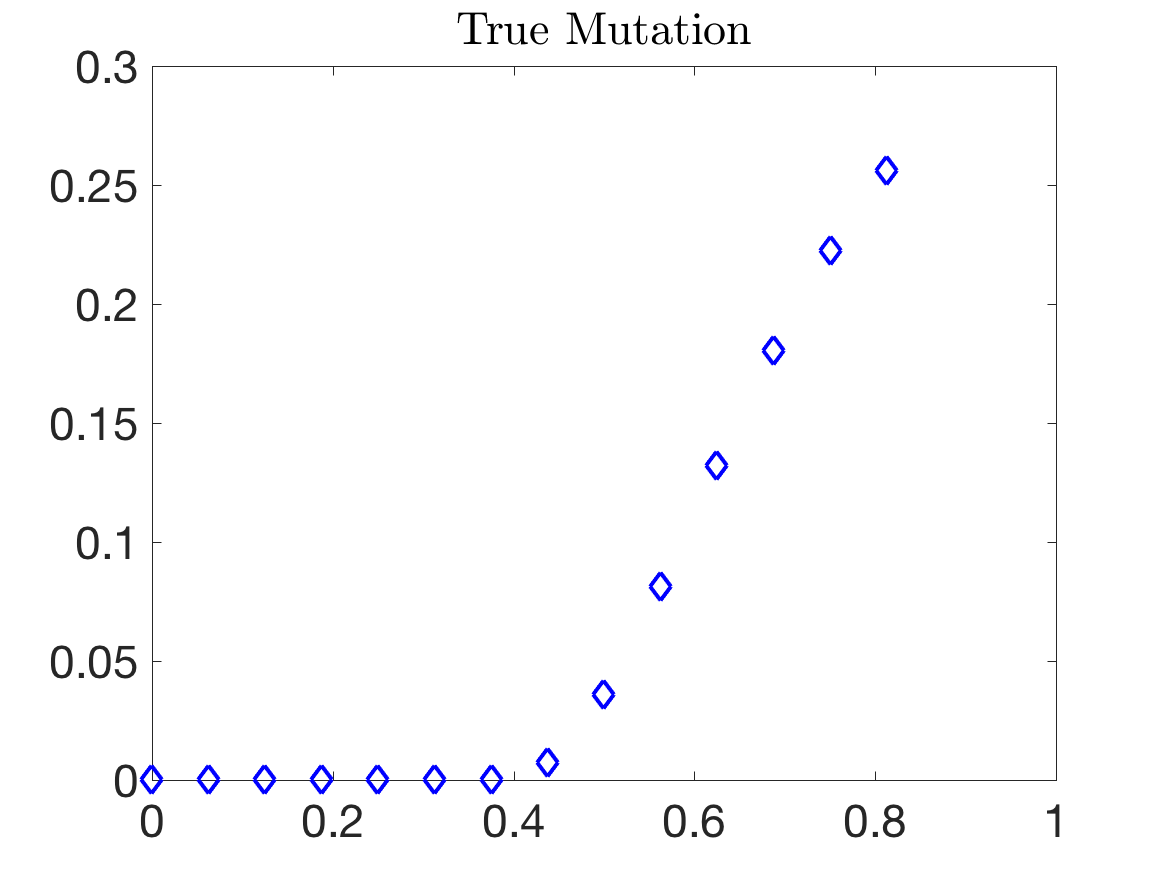}
	\includegraphics[width=0.45\textwidth]{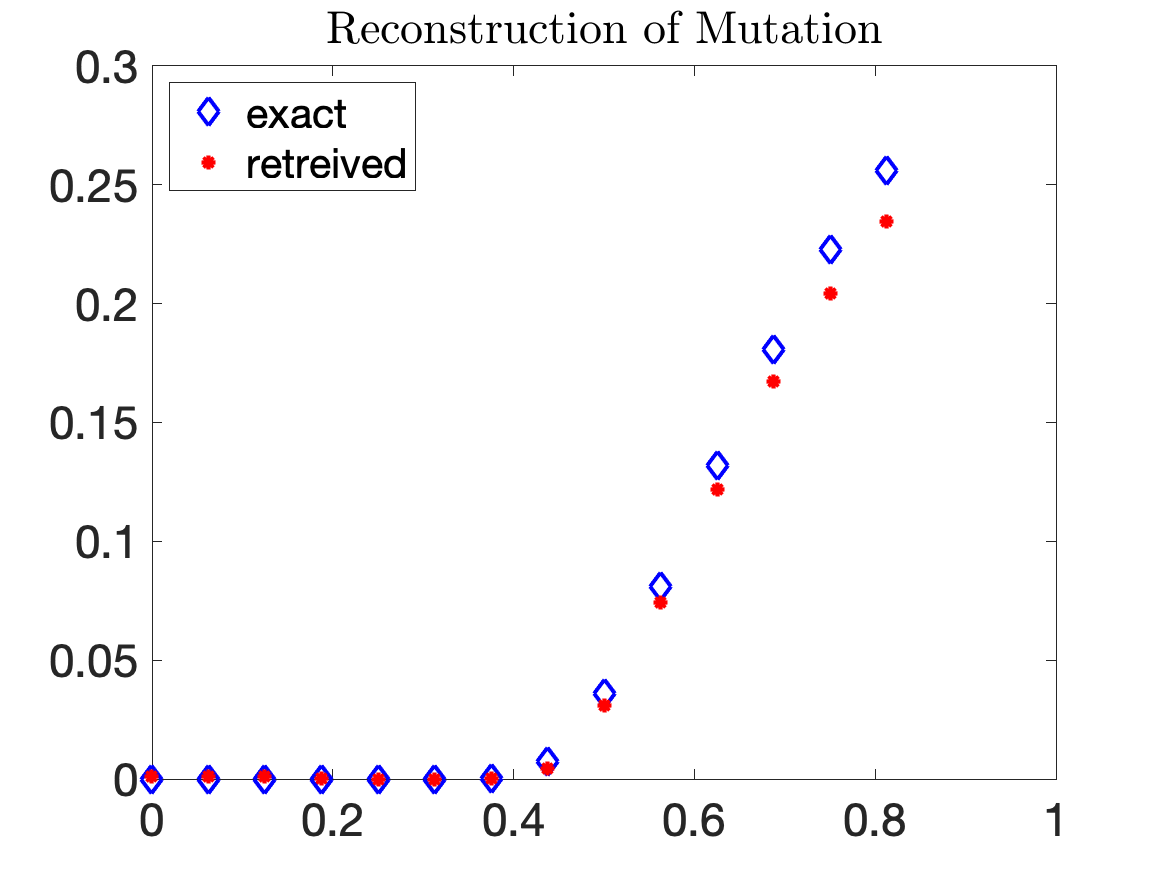}\\
	\vspace{-0.1cm}
         \includegraphics[width=0.45\textwidth]{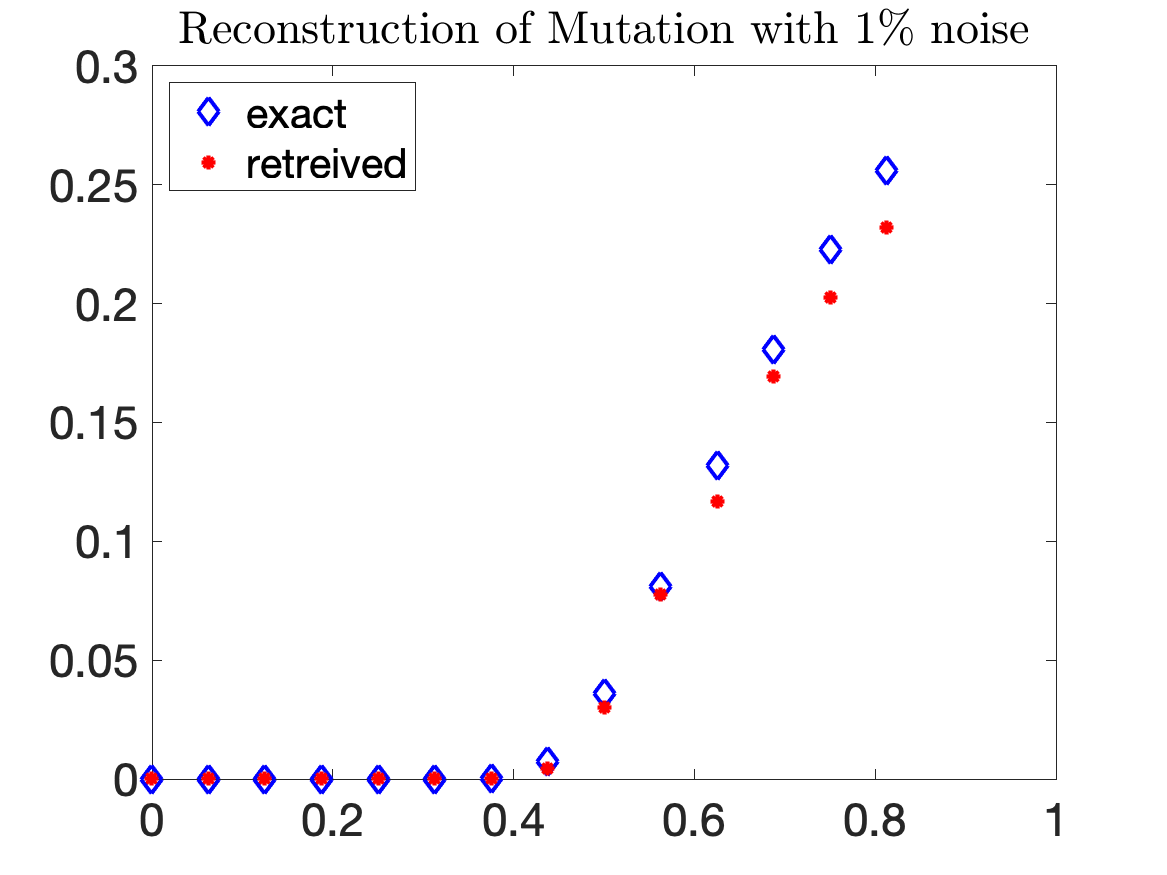}
         \includegraphics[width=0.45\textwidth]{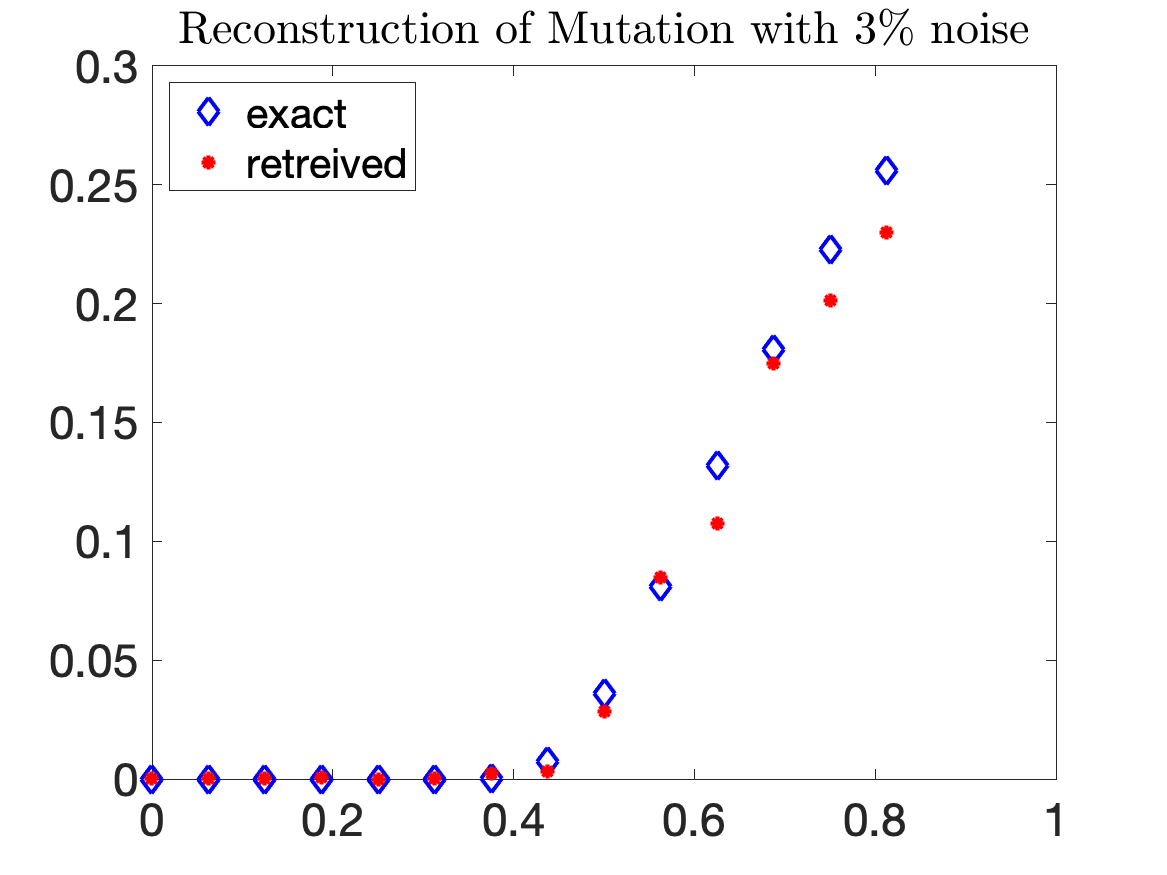}
	\caption{\small{\textit{Reconstruction of \dt{the unknown} mutation law \re{\dt{term $\widetilde{Q}_{2}(v)$ (that appears in mutation law given in \eqref{eq:Mutation_DepOn_V}} for} model \eqref{eq:Direct_Mutation_DepOn_C}. First row: \re{(left)}  the true mutation law, and  \re{(right)} the reconstructed mutation law restricted to $\A_{v}$ with no noise and $\alpha^{*}=10^{-12}$. Second row: \re{(left) the reconstructed} mutation law obtained for $1\%$ noisy data and $\alpha^{*}=10^{-12}$; and (right) \re{the reconstructed mutation law obtained for} $3\%$ noisy data and $\alpha^{*}=10^{-12}$. For all plots in this figure: 1) the first axis represents the the values for $v\in[\bar{v}^{min},\bar{v}^{max}]$; 2) second axis represents the values for mutation. \re{The numerical simulations are obtained using the parameters given in Table \ref{paramSetTable}. } }}} 
	\label{fig:True_Recon_Mut_DepOn_V}
\end{figure}

An acceptable numerical reconstruction of the mutation law $\dt{m}^{c_{1}^{*}, c_{2}^{*}, v^{*}}$, i.e., 
\bequ
\dt{m}^{c_{1}^{*}, c_{2}^{*}, v^{*}} := m^{\alpha^{*}},
\eequ  
is obtained for the choice of the regularisation parameter $\alpha^{*}$, which throughout this work is selected based on a standard discrepancy principle argument \cite{Morozov1984}.\\

Figure \ref{fig:True_Recon_Mut_Term_DepOn_C} shows the reconstruction of the cancer cell mutation law \dt{$\widetilde{\Q}_{1}(c_{1})$} \re{for} model \eqref{eq:Direct_Mutation_DepOn_C} in the presence of  the measurements given by \eqref{comp_meas_1} and \eqref{exactMeasurementsOnePopulation} that are considered here both exact and affected by a level of noise $\delta\in\{1\%, 3\%\}$. The first row of this figure shows from left to right the true mutation law restricted at the maximal accessible region $\A_{c}$ (where the reconstruction is being attempted) and the reconstruction of the mutation law equation \dt{$\widetilde{\Q}_{1}(c_{1})$}  \re{given by }\eqref{eq:Mutation_DepOn_C} on $\A_{c}$ with no noise, respectively. The second row of the figure show from left to right the reconstruction of the mutation law on $\A_{c}$ with $1\%$, and $3\%$ of noise in the measured data, respectively.

Figure \ref{fig:True_Recon_Mut_DepOn_V} shows the reconstruction of unknown nonlinear cancer cell mutation \dt{term }\re{\dt{$\widetilde{Q}_{2}(v)$ (that is part of cell mutation law introduced in \eqref{eq:Mutation_DepOn_V})} for model \eqref{eq:Direct_Mutation_DepOn_C}} in the presence of  the measurements given by \eqref{comp_meas_1} and  \eqref{exactMeasurementsOnePopulation} that are considered here both exact and affected by a level of noise $\delta\in\{1\%, 3\%\}$. Again, the first row shows from left to right the true mutation law restricted to $\A_{v}$ where the reconstruction is attempted, and the reconstruction of the mutation law equation \eqref{eq:Mutation_DepOn_V} on $\A_{v}$ with no noise, respectively. The second row of the figure show from left to right the mutation reconstruction on $\A_{v}$ with $1\%$, and $3\%$ of noise in the measured data, respectively.

From Figure \ref{fig:True_Recon_Mut_Term_DepOn_C} and Figure \ref{fig:True_Recon_Mut_DepOn_V} we observe that \re{when the measurement data is not affected by noise,} we obtain good mutation laws reconstructions in both cases, i.e., (i) mutation depends on primary tumour only, and (ii) mutation depend on the ECM only. However, as expected, as soon as the level of noise \re{in the measurements} increases, the reconstruction gradually looses accuracy in both cases (i.e., case (i) shown in  Figure \ref{fig:True_Recon_Mut_Term_DepOn_C} lower panels, and case (ii) shown in Figure \ref{fig:True_Recon_Mut_DepOn_V} lower panels). \re{This loss in accuracy increases with the increase in the $c_{1}$ density (horizontal axis in Figure \ref{fig:True_Recon_Mut_Term_DepOn_C}) and in the $v$ density (horizontal axis in Figure \ref{fig:True_Recon_Mut_DepOn_V}).}

\subsection{\re{Reconstruction of Nonlinear Mutation Law \dt{in case} (iii)} }
\label{Recon_Mutation_in_Local_Model}

In this section, we study the inverse problem of identifying the general mutation law \dt{in case } \re{(iii)}, namely $\Q(c_{1},v)$, \dt{which} is unknown \dt{and depends nonlinearly in an autonomous manner on both primary tumour cell density $c_{1}$ and on ECM density $v$}.
\subsubsection{Inverse Problem Setup: Forward Solver Computational Formulation}
\label{Computational_Approach_for_Direct_Problem}
\label{Inverse_Problem_Forward_Operator_and_Computational_Approach_CV}


\dt{Maintaining here the same spatio-temporal discretisation notations introduced in Section \ref{Inverse_Problem_Forward_Operator_and_Computational_Approach} for the cancer cell populations and ECM, we seek to identify the unknown nonlinear mutation law $\Q(c_{1},v)$ such that we match the measurements \eqref{eq:Additional_Information_Two_Pop}. To achieve this, similar to the previous two cases addressed in Sections \ref{Inverse_Problem_Forward_Operator_and_Computational_Approach}$-$\ref{Minimization_Algorithms}, within an suitable space of functions $\Mu^{2}$ that will be detailed below, we proceed to identify an appropriate a mutation approximating function defined on the 2-dimensional region where the pair $(c_{1}, v)$ ranges during its spatio-temporal evolution, namely $\overline{m}^{\,c_{1}^{*}, c_{2}^{*}, v^{*}}:[0,K_{c}]\times[0,K_{c}]\to [0,\infty)$. The function $\overline{m}^{\,c_{1}^{*}, c_{2}^{*}, v^{*}}$ will be identified such that this will directly determine an acceptable \emph{``mutation law candidate"}, denoted here by $\Q^{c_{1}^{*}, c_{2}^{*}, v^{*}}(\cdot,\cdot) $, which will enable a solution in \emph{forward model} \eqref{eq:Direct_Mutation_DepOn_C}  that matches the measurements \eqref{eq:Additional_Information_Two_Pop}. To select this mutation law candidate $\Q^{c_{1}^{*}, c_{2}^{*}, v^{*}}(\cdot,\cdot) $, adopting a similar approach to the one in cases $(i)$ and $(ii)$, we involve again an appropriately constructed \emph{mutation trial operator} $\M^{2}(\cdot,\cdot,\cdot):=\{\M^{2}_{i,j}(\cdot,\cdot,\cdot)\}_{i,j=1\dots N}$, $\M^{2}(\cdot,\cdot,\cdot):\R^{N\times N}\times \R^{N\times N}\times \Mu^{2}\to \R^{N\times N} $ that will be specified in a moment and will enable us the express the acceptable mutation law candidate as
%
\bequ
\label{ip3_sol}
\Q^{c_{1}^{*}, c_{2}^{*}, v^{*}}(\tilde c^{\overline{m}^{\,c_{1}^{*}, c_{2}^{*}, v^{*}}}_{1,i,j}(t),\tilde v^{\overline{m}^{\,c_{1}^{*}, c_{2}^{*}, v^{*}}}_{i,j}(t)):=\M^{2}_{i,j}(\tilde c^{\overline{m}^{\,c_{1}^{*}, c_{2}^{*}, v^{*}}}_{1}(t),\tilde v^{\overline{m}^{\,c_{1}^{*}, c_{2}^{*}, v^{*}}}(t),\overline{m}^{\,c_{1}^{*}, c_{2}^{*}, v^{*}})
\eequ 
where $\tilde c^{\overline{m}^{\,c_{1}^{*}, c_{2}^{*}, v^{*}}}_{1}(t):=\{\tilde c^{\overline{m}^{\,c_{1}^{*}, c_{2}^{*}, v^{*}}}_{1}(t)\}_{i,j=1\dots N}$ and $\tilde v^{\overline{m}^{\,c_{1}^{*}, c_{2}^{*}, v^{*}}}(t):=\{\tilde v^{\overline{m}^{\,c_{1}^{*}, c_{2}^{*}, v^{*}}}(t)\}_{i,j=1\dots N}$ represents the solution for the density of the primary cell population and of the ECM, respectively, which are obtained for model \eqref{eq:Direct_Mutation_DepOn_C} when, instead of the unknown term $\Q(\cdot,\cdot)$, in the mutation law we use the trial mutation term $\M^{2}(\cdot, \cdot, m^{c_{1}^{*},c_{2}^{*},v^{*}})$}

\dt{Denoting by $\G^{2}_{_{M}}:=\{(\eta_{l},\zeta_{k})\}_{l,k=1\dots M}$  the equally spaced grid given by the uniform discretisation of the $[0,K_{c}]\times[0,K_{c}]$ with step size $\Delta \eta=\Delta \zeta>0$, the space of functions where we seek to identify the function $\overline{m}^{\,c_{1}^{*}, c_{2}^{*}, v^{*}}$ is an $M\times M$- dimensional space of potential mutation law shape candidates, which is given by}
\dt{\bequ
\begin{split}
\Mu^{2}:=&\bigg\{\overline{m}:[0,K_{c}]\times[0,K_{c}]\to \R\,\bigg|\,\, \overline{m} \restrict{E_{l,k}}= \sum\limits_{p,q=0,1} \overline{m} (\eta_{l+p},\zeta_{k+q}) \phi_{l+p,k+q},\quad \forall E_{l,k}\in \G_{_{M}}^{2,\,tiles}\bigg\} \\
&\textrm{with} \quad \G_{_{M}}^{2,\,tiles}:=\{ E_{l,k}:=[\eta_{l},\eta_{l+1}]\times [\zeta_{k},\zeta_{k+1}]\,|\, l,k = 1\dots M-1\}, \,\, and\\
& \quad \qquad \forall \,\, E_{l,k}\in \G_{_{M}}^{2,\,tiles}, \,\,\{\phi_{l+p,k+q}\}_{p,q=0,1} \textrm{ are the usual bilinear shape functions on $E_{l,k}$.} 
\end{split}
\eequ}
\dt{Therefore, for any $\overline{m}\in \Mu^{2}$, the  \emph{trial proliferation operator} $\M^{2}$ is given by
\bequ
\label{attempted_prolif}
\begin{split}
&\qquad\qquad \M^{2}_{i,j}(\tilde c^{\overline{m}}_{1}(t),\tilde v^{\overline{m}}(t), \overline{m}):=  \overline{m}\restrict{E_{l,k}}(\tilde c^{\overline{m}}_{1}(t),\tilde v^{\overline{m}}(t))\\
& \textrm{with $(l,k)\in \Lambda_{i,j}:= \left\{l,\,k \!\in\! \{1,\dots, M\!-\!1\}\,\big| \, \exists E_{l'\!,k'}\!\in\G_{_{M}}^{tiles} \textrm{ such that }(\tilde c^{\overline{m}}_{1,i,j}(t),\tilde v^{\overline{m}}_{i,j}(t))\!\in\! E_{l'\!,k'} \!\right\}$,}\\
& \textrm{and noting also here that $(l,k)$ is independent of its choice within $\Lambda_{i,j}$}. 
\end{split}
\eequ}
Here, \dt{$\tilde c^{\overline{m}}_{1}(t):=\{\tilde c^{\overline{m}}_{1,i,j}(t)\}_{i,j=1\dots N}$, $\tilde c^{\overline{m}}_{2}(t):=\{\tilde c^{\overline{m}}_{2,i,j}(t)\}_{i,j=1\dots N}$ and $\tilde v^{\overline{m}}(t) := \{\tilde v^{\overline{m}}_{i,j}(t)\}_{i,j=1\dots N}$} represent the solutions at the grid points and time $t>0$ for the cancer cells and ECM densities obtained with model \dt{\eqref{eq:Direct_Mutation_DepOn_C}} when this uses $\M^{2}_{i,j}(\cdot,\cdot,\overline{m})$ as mutation law given in \eqref{attempted_prolif}. \dt{Therefore, in space-discretised form, model \eqref{eq:Direct_Mutation_DepOn_C}} can be recasted also in this case as 
\bequ
\label{eq:spatial_discretization}
\frac{\partial}{\partial{t}}
\begin{bmatrix}
\tilde c^{\overline{m}}_{1} \\
\tilde c^{\overline{m}}_{2} \\
\tilde v^{\overline{m}} 
\end{bmatrix}
= 
\left [
\begin{array}{l}
\H^{1}(\tilde c^{\overline{m}}_{1}, \tilde c^{\overline{m}}_{2}, \tilde{v}^{\overline{m}}, \overline{m})\\
\H^{2}(\tilde c^{\overline{m}}_{1}, \tilde c^{\overline{m}}_{2}, \tilde{v}^{\overline{m}}, \overline{m})\\
\H^{3}(\tilde c^{\overline{m}}_{1}, \tilde c^{\overline{m}}_{2}, \tilde{v}^{\overline{m}})\\
\end{array}
\right ], 
\eequ
\dt{Also in in this case (i.e., case $(iii)$) we have that $\H^{1}(\cdot,\cdot,\cdot,\cdot)$ and $\H^{2}(\cdot,\cdot,\cdot,\cdot$ are correspondingly defined through equations \eqref{opH1} and \eqref{opH1} when the \emph{trial mutation form} for the full mutation law (which in cases $(i)$ and $(ii)$ was given by operator $\overline{\M}_{r}$) is given here by $\M^{2}_{i,j}(\tilde c^{\overline{m}}_{1}(t),\tilde v^{\overline{m}}(t), \overline{m})$. } 
Finally, \dt{adopting again the same time discretisation as in Section \ref{Inverse_Problem_Forward_Operator_and_Computational_Approach} and using the Euler method, a time marching step} for system \eqref{eq:spatial_discretization} \dt{is given by following operator}
\bequ\label{MOL1}
\begin{split}
&\overline{K}_{\overline{m}}: \R^{N\times N}\times \R^{N\times N}\times \R^{N\times N} \to \R^{N\times N} \times \R^{N\times N} \times \R^{N\times N}\\
\textrm{given by}&\\
&\overline{K}_{\overline{m}}
\left(
\begin{bmatrix}
\tilde c^{m_{1}, n}_{1}\\[0.2cm] 
\tilde c^{m_{1}, n}_{2}\\[0.2cm]
\tilde v^{m_{1},n}
\end{bmatrix}
\right)
:= 
\begin{bmatrix}
\tilde c^{m_{1}, n}_{1}\\[0.2cm] 
\tilde c^{m_{1}, n}_{2}\\[0.2cm]
\tilde v^{m_{1},n}
\end{bmatrix}
+ \Delta{t} 
\begin{bmatrix}
\H^{1}(\tilde c^{\overline{m}, n}_{1}, \tilde c^{\overline{m}, n}_{2}, \tilde v^{\overline{m}, n}, \overline{m})\\[0.2cm]
\H^{2}(\tilde c^{\overline{m}, n}_{1}, \tilde c^{\overline{m}, n}_{2}, \tilde v^{\overline{m}, n}, \overline{m})\\[0.2cm]
\H^{3}(\tilde c^{\overline{m}, n}_{1}, \tilde c^{\overline{m}, n}_{2}, \tilde v^{\overline{m}, n})
\end{bmatrix},
\end{split}
\eequ
where, for any $n\in \{0,\dots, L\}$, we have \dt{$\tilde c^{\overline{m}, n}_{1} := \tilde c^{\overline{m}}_{1}(t_{n})$, $\tilde c^{\overline{m}, n}_{2} := \tilde c^{\overline{m}}_{2}(t_{n})$, and $\tilde{v}^{\overline{m},n}:=\tilde{v}^{\overline{m}}(t_{n})$}, \dt{while
\bequd  
\begin{array}{l}
\H^{1}(\tilde c^{\overline{m}, n}_{1}, \tilde c^{\overline{m}, n}_{2}, \tilde{v}^{\overline{m},n}, \overline{m}) := \H^{1}(\tilde c^{\overline{m}}_{1}(t_{n}), c^{\overline{m}}_{2}(t_{n}), \tilde{v}^{\overline{m}}(t_{n}), \overline{m}),\\[0.2cm]
 \H^{2}(\tilde c^{\overline{m}, n}_{1}, \tilde c^{\overline{m}, n}_{2}, \tilde{v}^{\overline{m},n}, \overline{m}) := \H^{2}(\tilde c^{\overline{m}}_{1}(t_{n}), c^{\overline{m}}_{2}(t_{n}), \tilde{v}^{\overline{m}}(t_{n}), \overline{m}),\\[0.2cm]
 \H^{3}(\tilde c^{\overline{m}, n}_{1}, \tilde c^{\overline{m}, n}_{2}, \tilde{v}^{\overline{m},n}):=\H^{3}(\tilde c^{\overline{m}}_{1}(t_{n}), \tilde c^{\overline{m}}_{2}(t_{n}), \tilde{v}^{\overline{m}}(t_{n})).
\end{array}
\eequd} 

\dt{Therefore, this allows us} formulate \emph{``forward operator"} \dt{denoted here by $\overline{K}$}  \dt{defined by
\bequ\label{MOL2}
\begin{split}
&\overline{K}:\Mu^{2}\to \R^{N\times N}\times \R^{N\times N}\times \R^{N\times N}\\
\textrm{given by}&\\
&\overline{K}(\overline{m}):=\underbrace{\overline{K}_{\overline{m}} \circ \overline{K}_{\overline{m}} \circ  \cdots \cdots  \circ \overline{K}_{\overline{m}} }_{L\, \text{times}}\left(
\begin{bmatrix}
\tilde{c}_{1,0}\\ 
\tilde{c}_{2,0}\\ 
\tilde{v}_{0}
\end{bmatrix}
\right)
\end{split}
\eequ
where $\tilde{c}_{1,0}$, $\tilde{c}_{2,0}$  and $\tilde{v}_{0}$ are the discretised initial conditions for \eqref{eq:Direct_Mutation_DepOn_C} (as introduced in \eqref{MOL_2}).} 

\dt{Similar to the previous two inverse problem cases, the \emph{forward operator} $\overline{K}$ will enable us to identify the $\overline{m}^{c_{1}^{*}, c_{2}^{*}, v^{*}}\in\Mu^{2}$ such that the solution of the resulting model matches measurements \eqref{eq:Additional_Information_Two_Pop}. Thus, similar to} Section \ref{Numerical_Solution_of_Inverse_Problem}, \dt{we observe again that} ${\overline{K}}$ \dt{can be written down as
\bequ
\label{map26APR2021}
\Mu^{2}\ni \overline{m}\longmapsto \overline{K}_{\overline{m}}\in \boldsymbol{\ell}^{2}(\boldsymbol{\ell}^{2}(\E\times \E\times \E); \boldsymbol{\ell}^{2}(\E\times \E\times \E))
\eequ}
Here,  $\boldsymbol{\ell}^{2}(\boldsymbol{\ell}^{2}(\E\times \E\times \E); \boldsymbol{\ell}^{2}(\E\times \E\times \E))$ \dt{is again} the finite-dimensional Bochner space of square integrable vector-value functions defined on $\boldsymbol{\ell}^{2}(\E\times \E\times \E)$ \dt{with} values in the same space.  \dt{As in Section \ref{Numerical_Solution_of_Inverse_Problem}, also here we have} that the mappings \dt{defined in} \eqref{map26APR2021} are continuous and compact. \dt{Therefore, } we obtain that $\overline{K}$ is also closed sequentially bounded, and \dt{as a consequence we have that also in this case we satisfy the} inverse problems hypotheses \dt{adopted} \cite{Engl_1989}, \dt{and thus, the convergence of the subsequent nonlinear Tikhonov regularisation is ensured. The Tikhonov functionals  $\{\overline{J}_{\alpha}\}_{ \alpha>0}$ will have essentially the same form as in the first two cases, except that the space where these are defined is different (i.e., in this case we have $\Mu^{2}$ of dimension $M\times M$ rather that $\Mu^{1}$ of dimension $M$ that we had in cases $(i)-(ii)$), namely:}
\bequ
\label{eq:Objective_Function_twoPop_2}
\begin{split}
&\overline{J}_{\alpha}: \Mu^{2} \to\R, \quad \forall \,\alpha>0,\\
\textrm{defined by}&\\
&\overline{J}_{\alpha}(\overline{m}) := \left\|\overline{K}(\overline{m})- \begin{bmatrix}
\tilde{c}_{1}^{*}\\ 
\tilde{c}_{2}^{*}\\
\tilde{v}^{*}
\end{bmatrix}\right\|^{2}_{2} + \alpha \|\,\overline{m}\,\|^{2}_{2} ,\quad \forall \,\overline{m}\in \Mu^{2}.
\end{split}
\eequ
\subsubsection{\re{Numerical Reconstruction of the Mutation Law $\Q(c_{1},v)$ \dt{in case} (iii) } }

\dt{Computationally, we address here the reconstruction of the \re{general mutation law $\Q(c_{1},v)$} that appears in case $(iii)$ for model \eqref{eq:Direct_Mutation_DepOn_C} in the presence of zero flux boundary conditions and the initial conditions prescribed in \eqref{eq:initial_condition_C1_1}-\eqref{eq:heterogenous_initial_condition_1}. Furthermore, the inverse problem is addressed in the presence of both exact and noisy measurements of the form detailed in \eqref{comp_meas_1}-\eqref{noise_in_the_measurements}, with the exact measurements given by} 
\bequ
 \label{exactMeasurementsTwoPopulation}
\tilde{c}^{*}_{1,exact}(x):= \bar{c}_{1}(x,t_{f}),\quad \tilde{c}^{*}_{2,exact}(x):= \bar{c}_{2}(x,t_{f}),\quad and \quad \tilde{v}^{*}_{exact}(x):=\bar{v}(x,t_{f}), \quad \forall x\in \Omega,
\eequ
\dt{where $\bar{c}_{1}(x,t), \,\, \bar{c}_{2}(x,t)$ and $\bar{v}(x,t)$ represent the solution densities for primary cell population, mutated cell population and ECM, respectively, which are obtained when model \eqref{eq:Direct_Mutation_DepOn_C} uses as mutation law the expression given by \eqref{eq:Mutation_DepOn_V} with the known term $\widetilde{\Q}_{2}(v)$ specified in \eqref{explicit_measurement_case_ii}}.
%
%

\dt{Although the dimensionality is different with respect to cases $(i)-(ii)$, the inversion method in this case follows the same steps and numerical minimisation procedure and steps that were described in Section \ref{Minimization_Algorithms} (reason for which we skip here details already outlined there). Indeed, for any $\alpha \in \{10^{-i}\,|\, i=1,\dots 12  \}$, the minimisation of the Tikhonov functional $J_\alpha$ \dt{starts with an initial guess} $\overline{m}_{0} = \mathbb{I} \times10^{-3}$, (where $\mathbb{I}$ represents the $M\times M$ matrix of ones), and leads to the numerical identification of the associated point of minimum $m^{\alpha}$ that correspond to the smalles mismatch between the associated solution (that is obtained  when model \eqref{eq:Direct_Mutation_DepOn_C} uses $m^{\alpha}$ as mutation law) and the measurements . Furthermore, since no data is available beyond the maximal region \emph{maximal accessible regions} $\A_{cv}$ given by} 
\bequd
\begin{array}{l}
\textrm{$\A_{cv}:=\re{\A_{c}\times \A_{v}=} [\bar{c}_{1}^{min},\bar{c}_{1}^{max}]\times [\bar{v}^{min},\bar{v}^{max}]$, \re{with $\A_{c}$ and $A_{v}$ described in \eqref{EqAc}-\eqref{EqAv}.} }\\[0.2cm]
\end{array}
\eequd
\dt{the reconstruction of the mutation law is explored only on $\A_{cv}$. Finally, using again an discrepancy principle-based argument to choose the regularisation parameter $\alpha*\in \{10^{-i}\,|\, i=1,\dots 12  \}$, the reconstructed mutation law will be given by the corresponding $M\times M$ matrix $\overline{m}^{\,c_{1}^{*}, c_{2}^{*}, v^{*}}:=\overline{m}^{\alpha*}$, which in turn will determine the precise shape of the reconstructed mutation law $\Q^{c_{1}^{*}, c_{2}^{*}, v^{*}}(\tilde c^{\overline{m}^{\,c_{1}^{*}, c_{2}^{*}, v^{*}}}_{1,i,j}(t),\tilde v^{\overline{m}^{\,c_{1}^{*}, c_{2}^{*}, v^{*}}}_{i,j}(t))$ that is defined as per equation \eqref{ip3_sol}.}
%

\begin{figure}[!ht]
	\centering
        \setlength{\arrayrulewidth}{0.5pt}
        \begin{tabular}{lc}
        a) & \\
            &\hspace{-0.6cm}
	\includegraphics[width=0.65\textwidth]{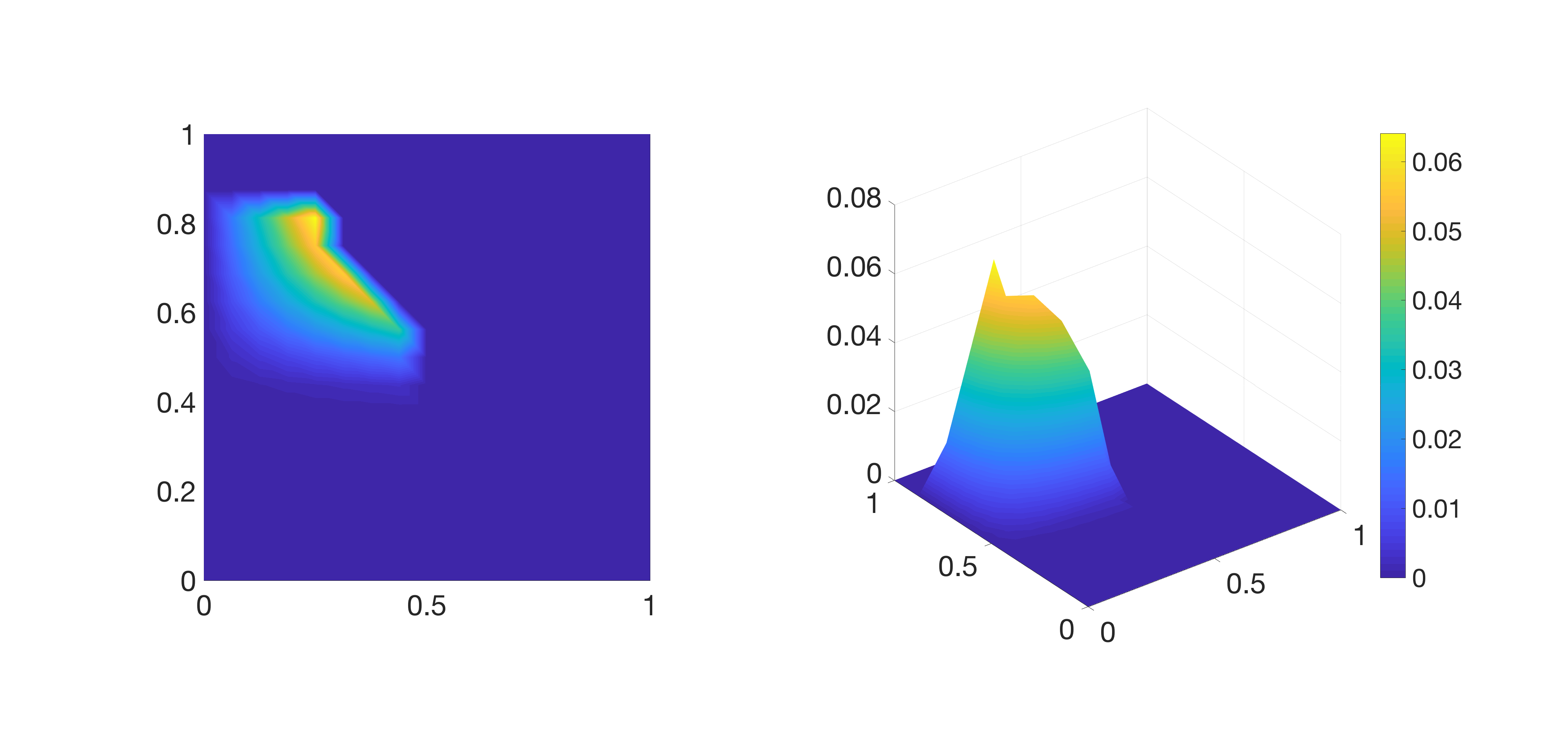}\\[-0.2cm]
         b)&\\
         &\hspace{-0.5cm}
	\includegraphics[width=0.3\textwidth]{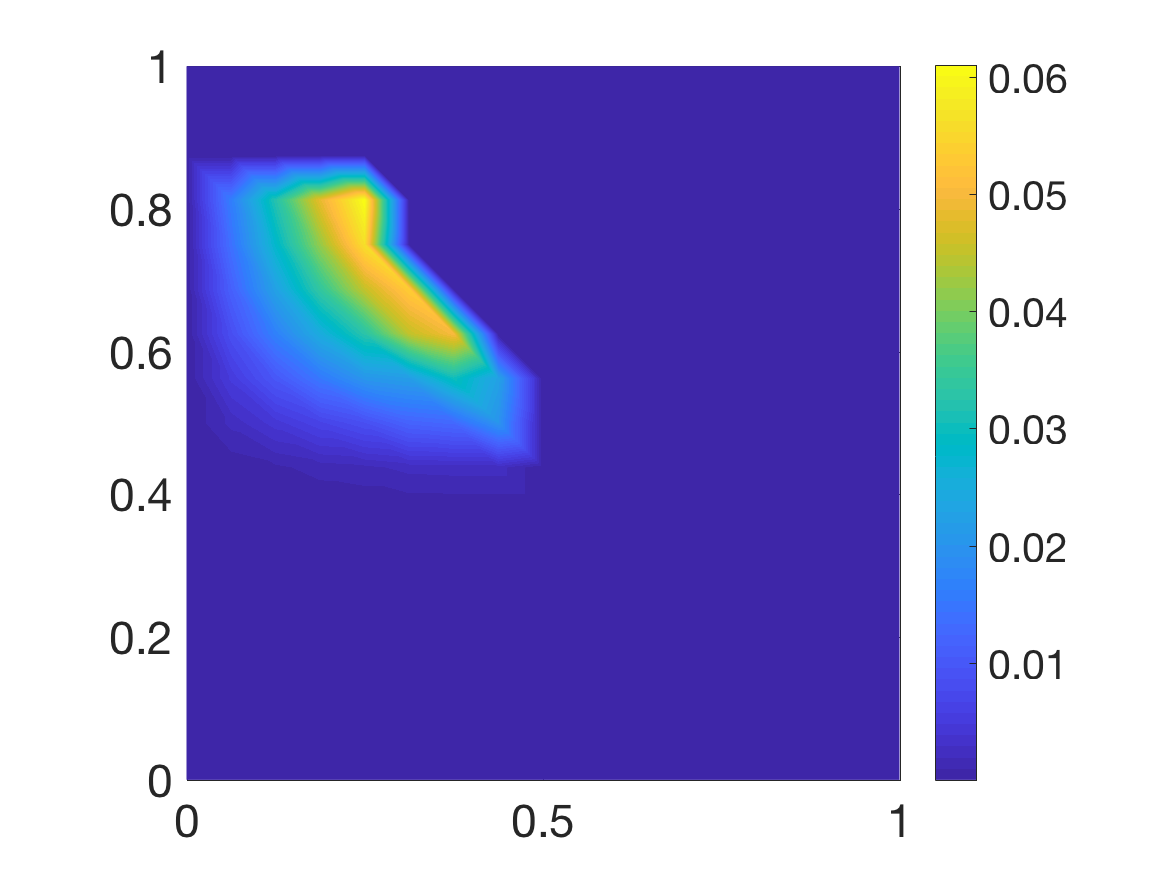}\hspace{-0.5cm}
	\includegraphics[width=0.3\textwidth]{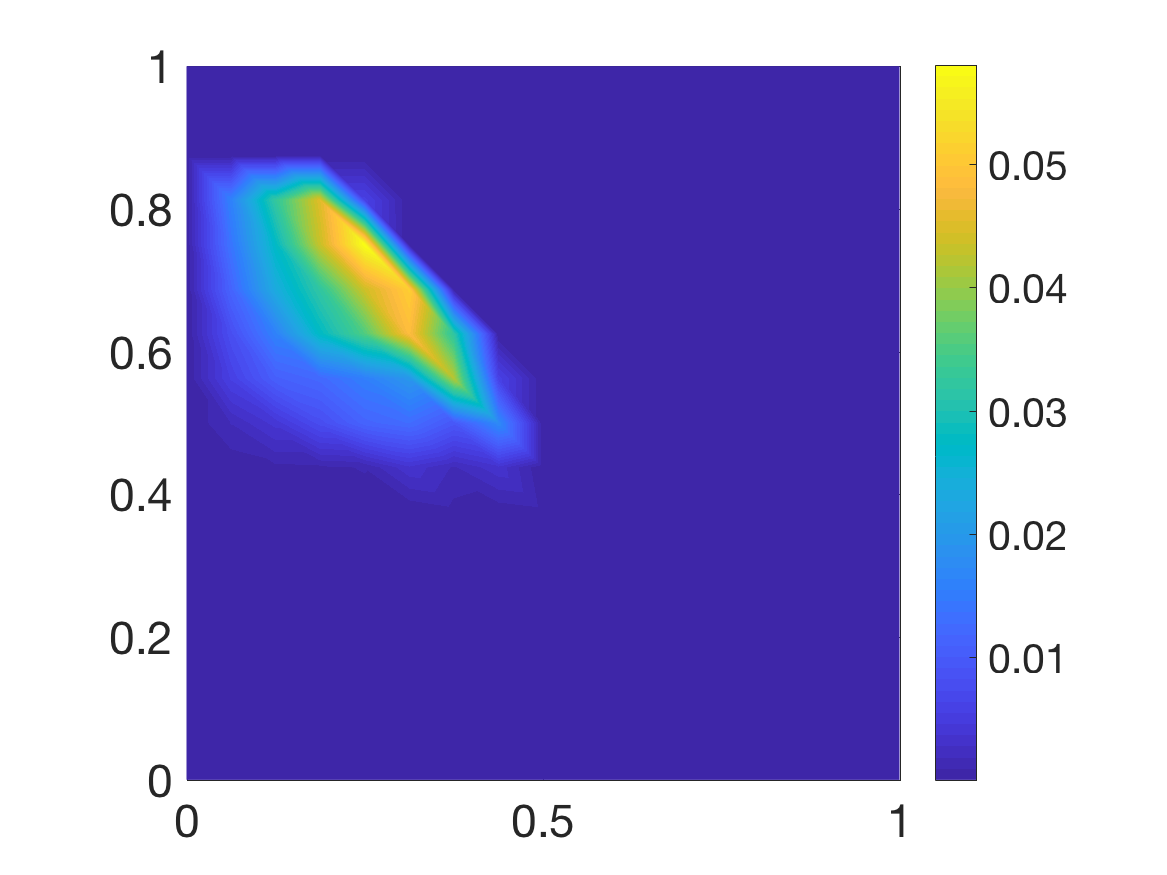}\hspace{-0.5cm}
	\includegraphics[width=0.3\textwidth]{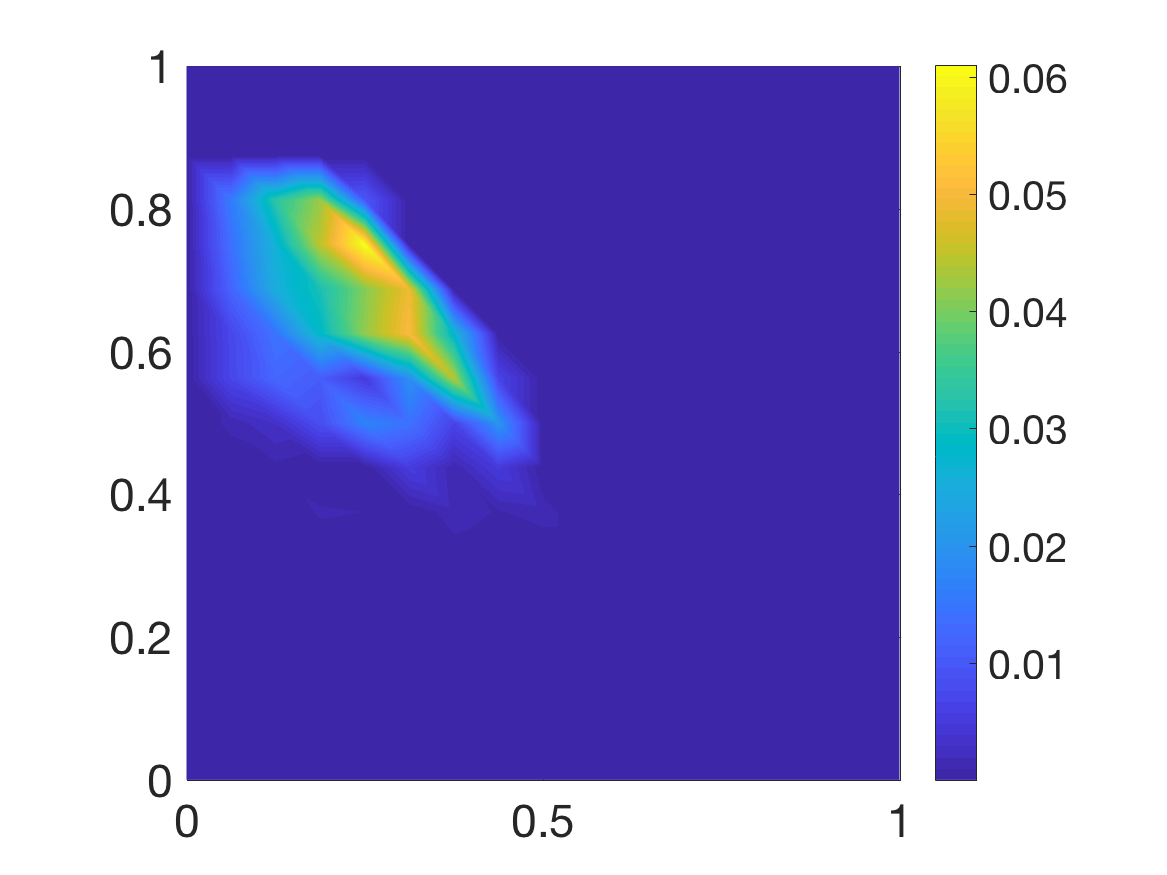}
	\end{tabular}
	\caption{\small{\textit{Reconstruction of \re{the general mutation law $\Q(c_{1},v)$ \dt{in case $(iii)$ }for} model \eqref{eq:Direct_Mutation_DepOn_C}. Row a) the true mutation law restricted to $\A_{cv}$; Row b) the reconstructed mutation law on $\A_{cv}$ in the presence of exact and noisy data: (left) exact data and $\alpha^{*}=10^{-8}$; (centre) $1\%$ noisy data and $\alpha^{*}=10^{-5}$; and (right) $3\%$ noisy data and $\alpha^{*}=10^{-5}$. For all plots in this figure: 1) first axis \re{shows} the values $c_{1}\in[\bar{c}_{1}^{min}, \bar{c}_{1}^{max}]$; 2) second axis \re{shows} the values $v\in [\bar{v}^{min},\bar{v}^{max}]$; and 3) colour bars represent the magnitude of mutation law or its reconstructions at each $(c_{1}, v)\in \A_{cv}$. \re{Numerical simulations are obtained using the parameters given in Table \ref{paramSetTable}. }}}} 
\label{fig:Recon_of_Mut_DepOnCV_shap_4}
\end{figure}

Figure \ref{fig:Recon_of_Mut_DepOnCV_shap_4} shows the reconstruction of the \re{most general cancer cell  mutation law $\Q(c_{1},v)$, starting from} the measurements given by \eqref{comp_meas_1} and  \eqref{exactMeasurementsOnePopulation} \re{for exact data as well as data impacted by} noise $\delta\in\{1\%, 3\%\}$. For comparison, the first row of this figure shows the true mutation law restricted at the maximal accessible region  $\A_{cv}$ where the reconstruction is being attempted. The second row of the figure show from left to right the reconstruction of the mutation law on $\A_{cv}$ with no noise \re{(left)}, with $1\%$ \re{noise (center)}, and $3\%$ noise \re{(right)} in the measured data. From this figure, we observe that we obtain a good \re{reconstruction of the mutation law} when the measurement data are not affected by noise. However, as expected, as soon as the level of noise increases in the measurements, the reconstruction gradually looses accuracy. 
%
%
%
%
%
\subsubsection{Reconstruction of Unknown Mutation Law \re{$\Q(c_{1},v)$ in case $(iii)$ for a Different Cell} Proliferation Rule for $c_{2}$}
\label{SectLocalMixedGrowth}

\re{Throughout the previous sections we assumed that both cancer cell populations proliferate logistically: $\mu_{c} c_{1,2}\left(1\frac{c_{1}+c_{2}+v}{K_{c}}\right)$. However, the sigmoid shape of tumour growth that is given by the logistic term can be obtained also with other proliferation rules, such as the Gompertz rule~\cite{Guiot2003}: $\mu_{c}c_{1,2}\log \left(\frac{K_{c}}{c_{1}+c_{2}+v} \right)$. This raises the question as to what happens when different cancer cell populations use different proliferation laws.}\\

\begin{figure}[!ht]
\centering
\setlength{\arrayrulewidth}{0.5pt}
\begin{tabular}{lc}
a) &\\
&\hspace{-0.6cm}
\includegraphics[width=0.65\textwidth]{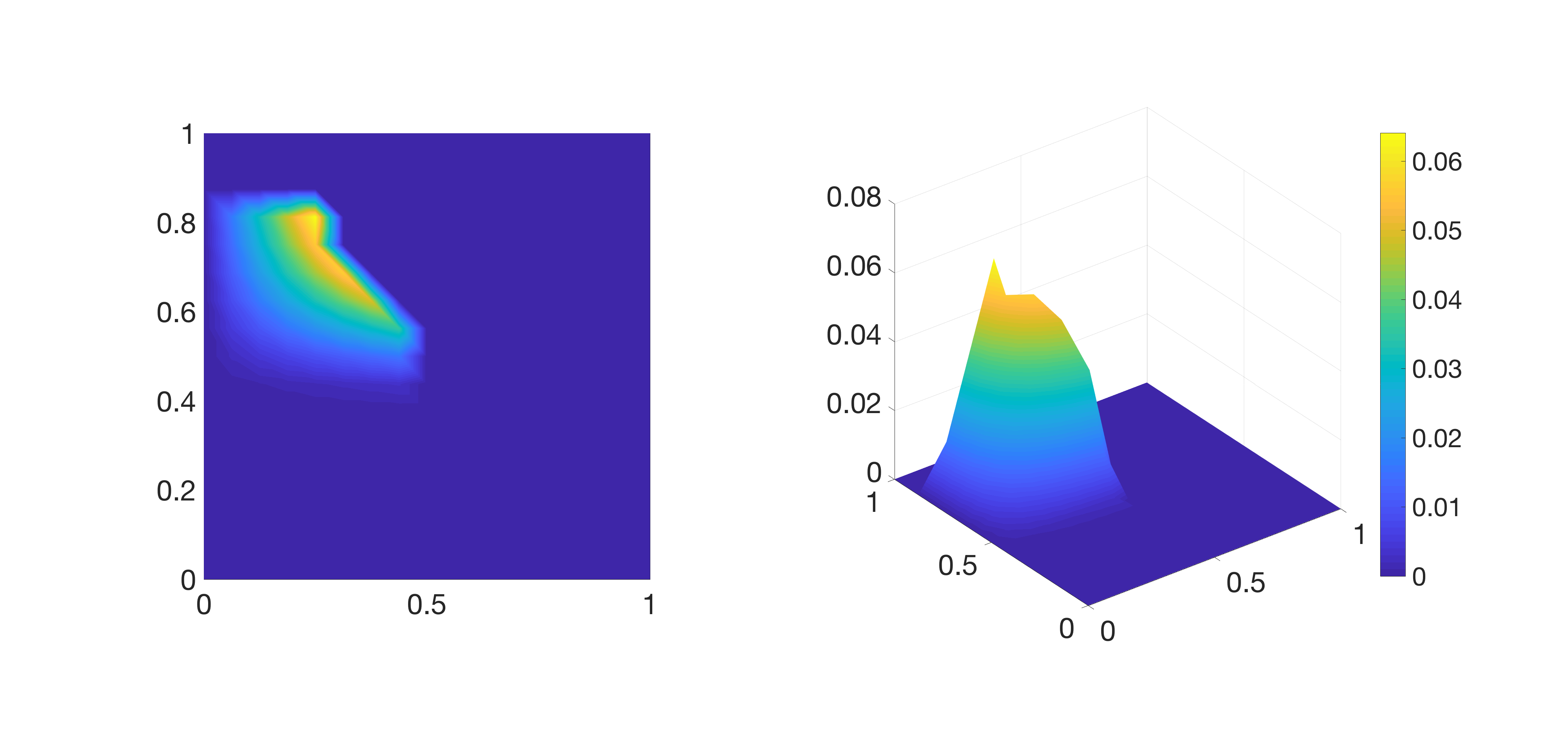}\\ 
b)& \\
&\hspace{-0.5cm}
\includegraphics[width=0.3\textwidth]{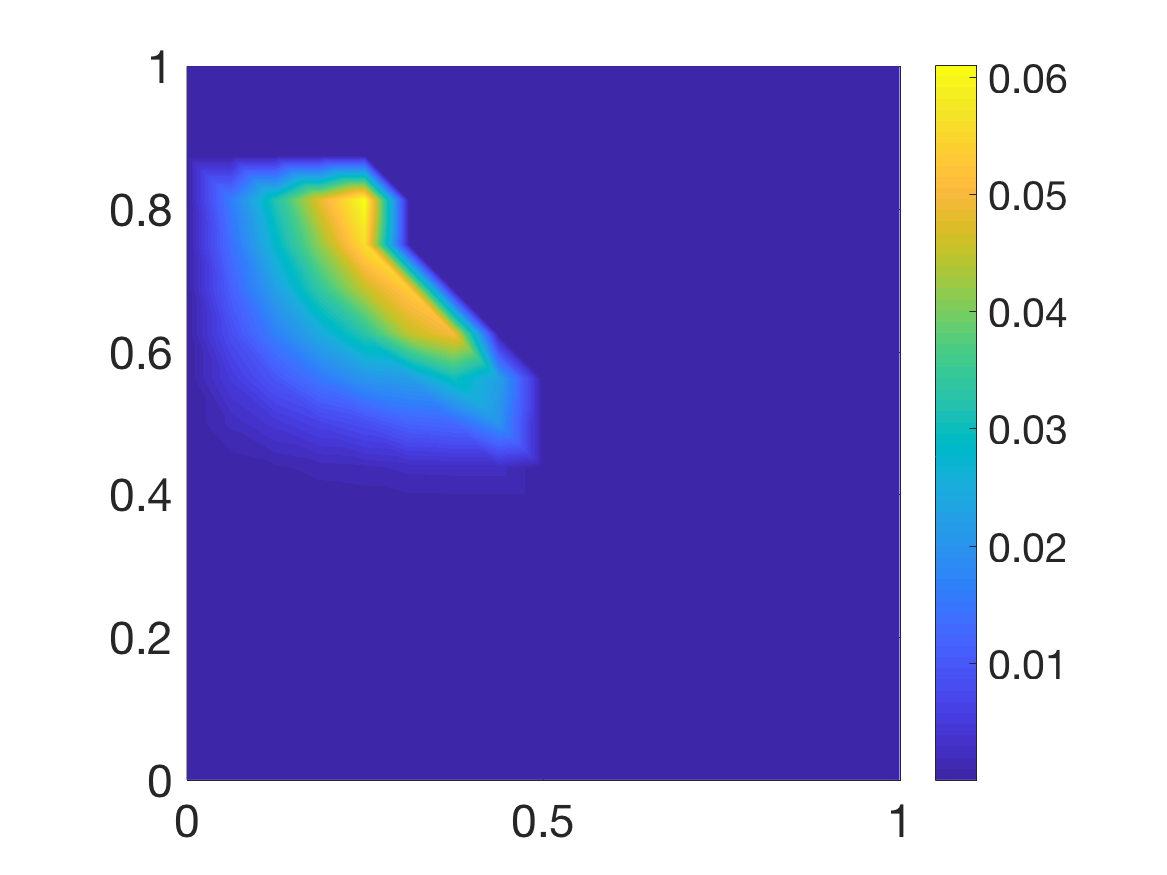}\hspace{-0.3cm}
\includegraphics[width=0.3\textwidth]{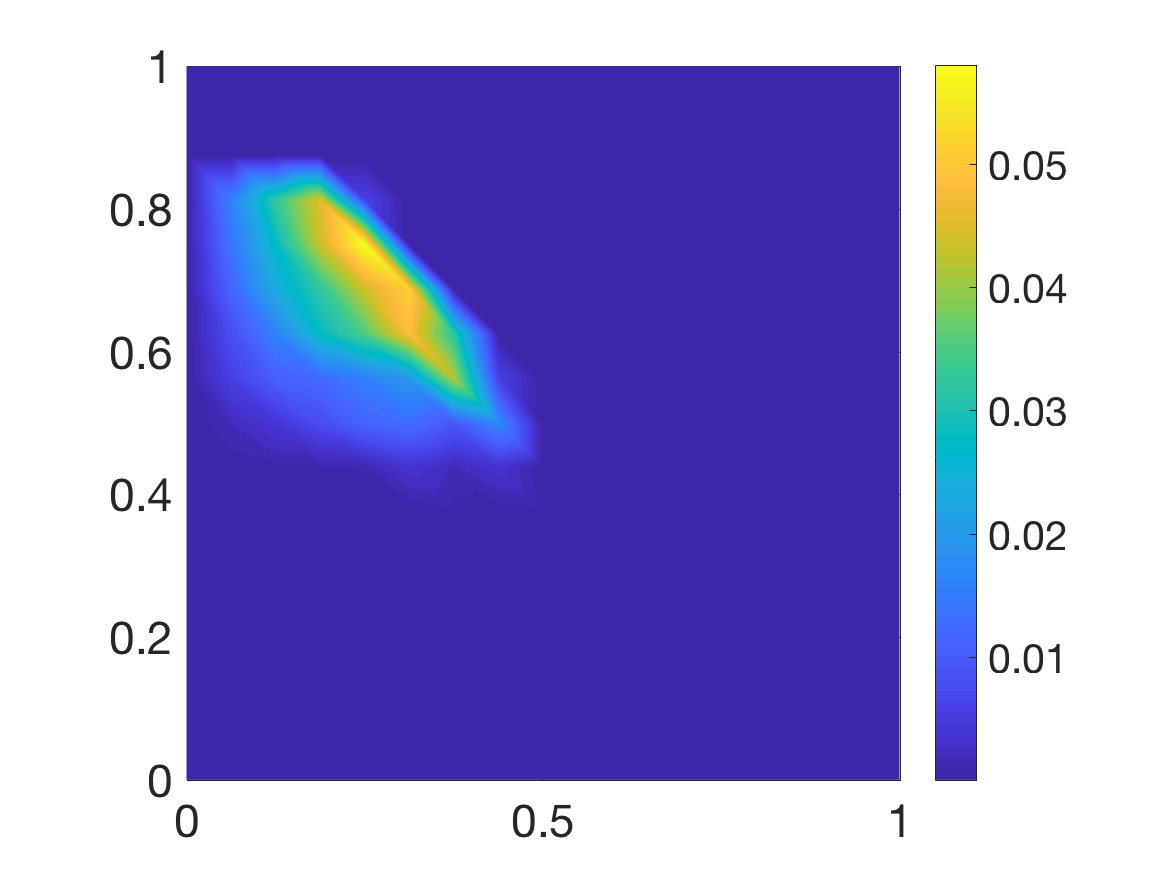}\hspace{-0.3cm}
\includegraphics[width=0.3\textwidth]{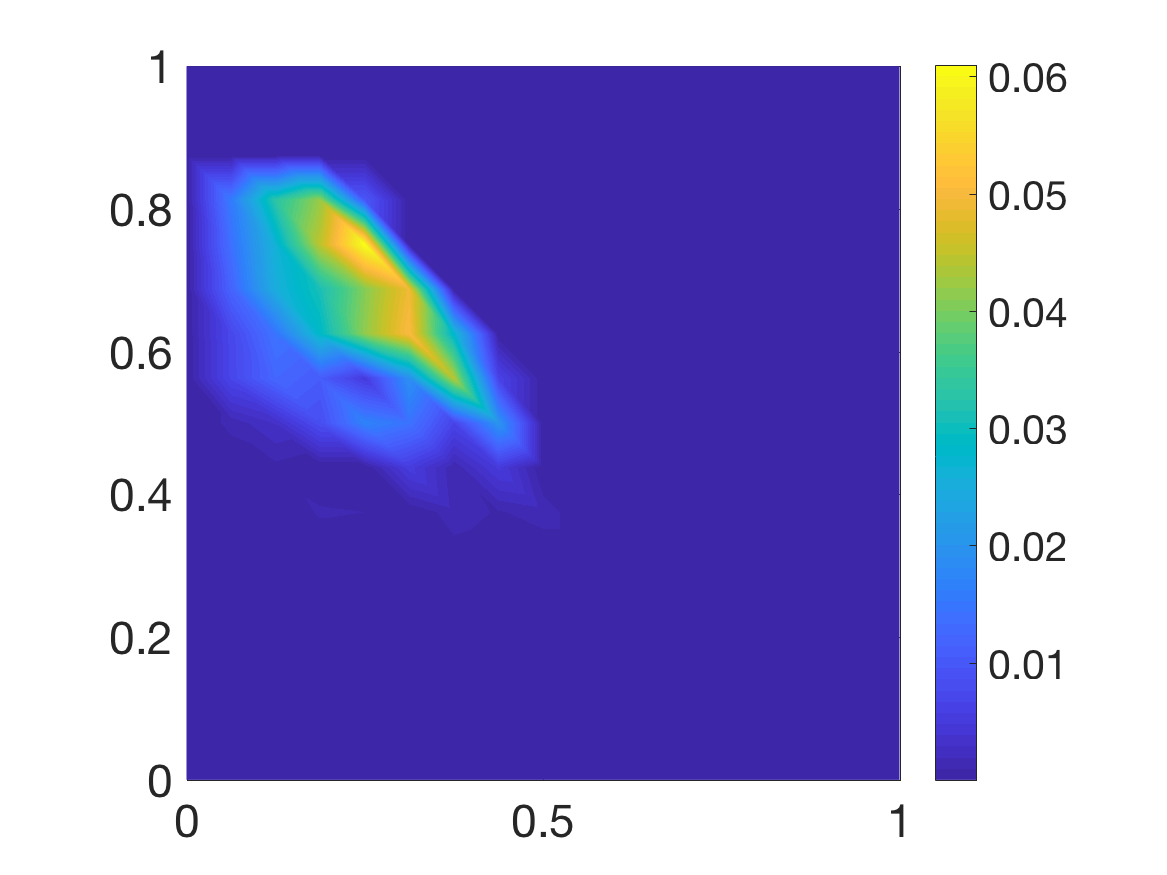}
\end{tabular}
\caption{Reconstruction of \re{general mutation law $\Q(c_{1},v)$ \dt{in case $(iii)$ }for} model \eqref{eq:Direct_Mutation_DepOn_C}, obtained \re{when using logistic cell proliferation for $c_{1}$ and Gompertz cell proliferation for $c_{2}$.} Row a): the true mutation law restricted to $\A_{cv}$; Row b) the reconstructed  mutation law on $\A_{cv}$ in the presence of exact and noisy data: (left) exact data and $\alpha^{*}=10^{-4}$; (centre) $1\%$ noisy data and $\alpha^{*}=10^{-4}$; and (right) $3\%$ noisy data and $\alpha^{*}=10^{-4}$. For all plots in this figure: 1) first axis \re{shows} the values $c_{1}\in[\bar{c}_{1}^{min}, \bar{c}_{1}^{max}]$; 2) second axis \re{shows} the values $v\in [\bar{v}^{min},\bar{v}^{max}]$; and 3) the colour bars represent the magnitude of mutation law or its reconstructions at each $(c_{1}, v)\in \A_{cv}$. \re{The numerical simulations are obtained with the parameters given in Table \ref{paramSetTable}. }}
\label{fig:Reconstruction_of_Mutation_Law_Mixed_Prol}
\end{figure}

In this subsection we reconstruct the general mutation law \re{$\Q(c_{1},v,t)$ \dt{in case $(iii)$} when we assume that the primary $c_{1}$ population  proliferation proliferates logistically, while the secondary $c_{2}$ population proliferates according to Gompertz law. In Figure \ref{fig:Reconstruction_of_Mutation_Law_Mixed_Prol} we present the numerical reconstruction results. We observe that the results are similar to those in Figure \ref{fig:Recon_of_Mut_DepOnCV_shap_4}; this could be explained by the fact that the mutation starts in the primary tumour which proliferate logistically.}

%
%
%
%
%
%
%
%
\section{\re{A Tumour Invasion Model with Nonlocal Dynamics}}\label{Sect6}

\re{Since various mathematical studies have assumed nonlocal cell-cell and cell-ECM interactions to explain the invasion of cancer cells~\cite{Armstrong_et_al_2006,Dumitru_et_al_2013,Domschke_et_al_2014,Gerisch_Chaplain_2008}, in the following we generalise the model \eqref{eq:Direct_Mutation_DepOn_C} by replacing the local haptotaxis towards local ECM gradients with a nonlocal haptotaxis flux generated by these cell-cell and cell-ECM interactions. }

\re{As before, we consider a primary cancer cell subpopulation $c_{1}(x,t)$ and a secondary mutated cancer cell subpopulation $c_{2}(x,t)$. These cancer cell populations interact with each other as well as with the extracellular matrix (ECM), $v(x,t)$, which they degrade and remodel.}  For compact notation, we consider the combined vector of primary cancer cells $c_{1}$, mutated cancer cells $c_{2}$ and extracellular matrix $v$ defined as
\[
\bu(x,t) := \left[c_{1}(x,t), c_{2}(x,t), v(x,t)\right]^{T}. 
\]
\re{We use this vector to describe in a compact manner the flux term generated by the nonlocal cell-cell and cell-ECM interactions ($\A_{1,2}(x,t,\bu{\cdot,t)}$)}. \dt{Therefore the coupled tumour dynamics in this case is given by:}
%
\begin{subequations} 
\label{eq:TwoPop_Nonlocal_Chap5}
\begin{align}
\label{eq:IP_Twopop_NonLocal_C1_Chap5}
\frac{\partial{c_{1}}} {\partial {t}} & = \underbrace{D_{1} \Delta c_{1} }_{\text{random motility}} - \underbrace{\nabla \cdot [c_{1}\mathcal{A}_{1}\left(x,t,\bu\left(\cdot,t\right)\right)]}_{\text{adhesion}} + \underbrace{\mu_{c} {c_{1}} \left( 1- \frac{c_{1}+c_{2} + v}{K_{c}} \right)}_{\text{logistic proliferation}} - \dt{\underbrace{\omega(t)}_\text{\parbox[t]{0.4in}{mutation\\[-0.2cm] switch}}\underbrace{\Q(\cdot,\cdot)}_\text{\parbox[t]{0.4in}{unknown\\[-0.2cm] mutation}}},\\
\label{eq:IP_Twopop_NonLocal_C2_Chap5}
\frac{\partial{c_{2}}} {\partial {t}} & = \underbrace{D_{2} \Delta c_{2} }_{\text{random motility}} - \underbrace{\nabla \cdot [c_{2}\mathcal{A}_{2}\left(x,t,\bu\left(\cdot,t\right)\right)]}_{\text{adhesion}} + \underbrace{\mu_{c} {c_{2}} \left( 1 - \frac{c_{1}+c_{2} +v}{K_{c}} \right)}_{\text{logistic proliferation}} + \dt{\underbrace{\omega(t)}_\text{\parbox[t]{0.4in}{mutation\\[-0.2cm] switch}}\underbrace{\Q(\cdot,\cdot)}_\text{\parbox[t]{0.4in}{unknown\\[-0.2cm] mutation}}},\\
\label{eq:IP_Twopop_NonLocal_V_Chap5}	
\frac{\partial {v}} {\partial {t} } & = \underbrace{- \rho (c_{1}+c_{2}) v}_{\text{degradation}} + \underbrace{\mu_{v} \left(K_{c} - c_{1}-c_{2} - v \right)^{+}}_{\text{ECM remodelling}}.
\end{align}
\end{subequations}
\re{As before, $D_{1,2}$ are the diffusion coefficients, $\mu_{c}$ is the cancer cell proliferation rate, $\rho$ is the ECM degradation rate, and $\mu_{v}$ is the ECM remodelling rate.} The flux term $\mathcal{A}_{p}\left(x,t,\bu\left(\cdot,t\right)\right)$, $p=1,2$, has been proposed in \cite{Domschke_et_al_2014,Gerisch_Chaplain_2008} to describe the directed movement of cells due to cell-cell and cell-matrix adhesion:
\begin{equation}
\mathcal{A}_{p}\left(x,t,\bu\left(\cdot,t\right)\right) := \frac{1}{R}\int\limits_{\textbf{B}\left(\left(0,0\right),R\right)} n\left(y\right) \cdot \mathcal{K}\left(\|y\|_{2}\right) \cdot g_{\dt{p}}\left(\bu\left(x+y,t\right),t\right)\chi_{_{\Omega}}\left(x+y\right) \ dy.
\end{equation}
\re{It is assumed that the interactions between a cell and its neighbouring cells as well as the components of the ECM occur inside a sensing region $\textbf{B}\left(\left(0,0\right),R\right) \subset \mathbb{R}$, where $R>0$ is the \textit{sensing radius}. }
In the above equation, $\chi_{_{\Omega}}\left(\cdot\right)$ represents the characteristic function of $\Omega$. Further, $n\left(y\right)$ is the \re{unit radial vector} given by 
\begin{equation} 
n\left(y\right):= 
\left\{
\begin{array}{lr}
y / \|y\|_{2} \quad \quad \quad  \text{if } \quad y \in \textbf{B}\left(0,R\right) \setminus \{\left(0,0\right)\},\\
\left(0,0\right) \quad \quad \text{otherwise}.
\end{array}
\right.
\end{equation}
\dt{Furthermore, $g_{\dt{p}}\left(\bu(x+y,t\right),t)$ represents the $p-$\emph{th} component of the adhesion function} $g\left(\bu\left(x+y,t\right),t\right)$ that is given by 
\begin{equation}
\begin{array}{c}
g\left(\bu,t\right) = [\textbf{S}_{cc} \bc + \textbf{S}_{cv} v] \cdot \left(K_{c} - c_{1}-c_{2} - v\right)^{+},\\[0.5cm]
\textrm{with}\quad 
\textbf{S}_{cc} =
\begin{pmatrix}
S_{c_{1}c_{1}} & S_{c_{1}c_{2}} \\
S_{c_{2}c_{2}}    & S_{c_{2}c_{2}}
\end{pmatrix}
\quad\quad 
\text{and}
\quad\quad 
\textbf{S}_{cv} =
\begin{pmatrix}
S_{c_{1}v} & 0 \\
0    & S_{c_{2}v}
\end{pmatrix}.
\end{array}
\end{equation}
Here $\{S_{c_{i}c_{j}}\}_{i=1,2}$ are the non-negative cell-cell adhesion \dt{strengths of the adhesion bonds established between the primary and mutated cancer cell populations, while  $\{S_{c_{i}v}\}_{i=1,2}$} are the non-negative cell-matrix adhesion \dt{stands for the strengths for the adhesive interactions between each of the two cancer subpopulation and the ECM}. \dt{Finally, no overcrowding} of the cell population and ECM \dt{over the tumour domain is ensured here through the term} $(K_{c} - c_{1}-c_{2} - v)^{+}:= \max(K_{c} - c_{1}-c_{2} - v, 0)$. 

The general mutation law $\Q(c_{1}, v)$ is assumed to be unknown and its identification will be our main focus in this section. Finally, \re{the dynamics described by} \eqref{eq:IP_Twopop_NonLocal_C1_Chap5}-\eqref{eq:IP_Twopop_NonLocal_V_Chap5} takes place in the presence of initial conditions \eqref{subeq:Initial_Conditions_One_Pop_1} and zero Neumann boundary conditions \eqref{eq:DP_BC}.
%
%

%
%
\subsection{Numerical Approach for the \re{Forward} Nonlocal Cancer Invasion Model}
\label{Numerical_Approach_for_Direct_Non_Local_Cancer_Model}
\begin{figure}[t]
	\centering
	\includegraphics[width=0.75\textwidth]{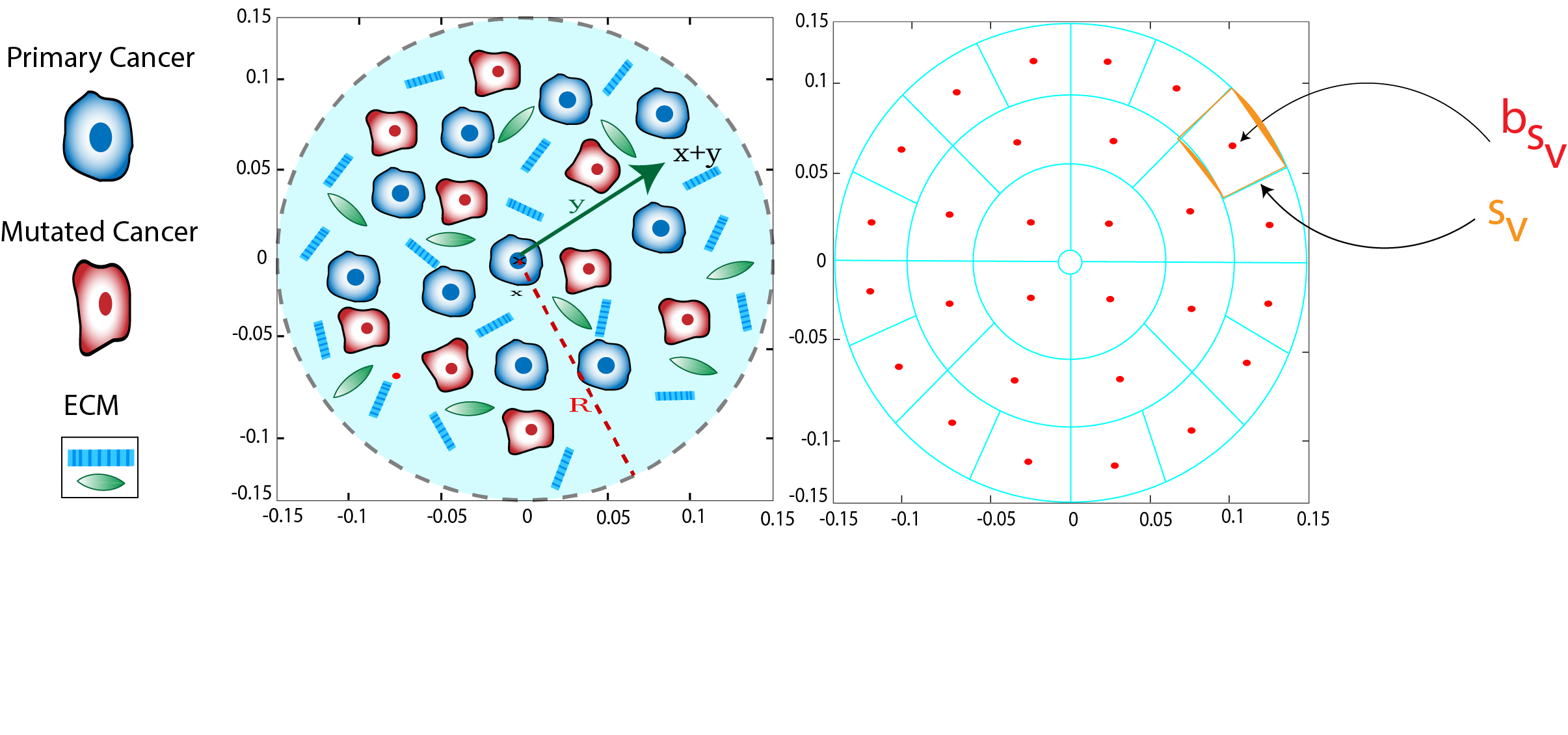}\\[-1.3cm]
	\caption{\small{\textit{Schematic shows the sensing region $B(0,R)$. Left figure describes the process of the cancer cells adhesion, the ball $B(x,R)$ centred at x and of radius R, the point $x+y$ with the direction vector $n(y)$ in green. Right figure shows the decomposition of the region using annulus sectors $S_{\nu}$ with barycentres $b_{S_{\nu}}$, highlighted with red dots. }}}
\label{fig:Schematic_One}
\end{figure}

In this section, we briefly discuss the numerical methods used to solve the \re{forward} model  \eqref{eq:IP_Twopop_NonLocal_C1_Chap5}-\eqref{eq:IP_Twopop_NonLocal_V_Chap5}. \re{To discretise the system in space we use the method of lines (MOL) approach}. \re{For the time-evolution of the system}, we use a  predictor-corrector scheme introduced in \cite{Dumitru_et_al_2013}, using the Euler method as the predictor and a trapezoidal-type rule as the corrector. In the reaction-diffusion equations \eqref{eq:IP_Twopop_NonLocal_C1_Chap5} and \eqref{eq:IP_Twopop_NonLocal_C2_Chap5}, the terms \re{$D_{p} \Delta c_{p}- \nabla \cdot [c_{p} \mathcal{A}_{p}\left(x,t,\bu\left(\cdot,t\right)\right)]$, $p=1,2$,} will be approximated through a second-order mid-point rule \cite{Leveque_2007_finite} as detailed below. In brief, \re{$\forall p\in\{1,2\}, \; n=\overline{0,L}, \; i,j=\overline{1,N}$,} these midpoint approximations are given by: 
\begin{equation}
\begin{cases}
c^n_{p,i,j+\frac{1}{2}}:=\frac{c^n_{p,i,j} + c^n_{p,i,j+1}}{2},\\ 
c^n_{p,i,j-\frac{1}{2}}:=\frac{c^n_{p,i,j} + c^n_{p,i,j-1}}{2},\\
c^n_{p,i+\frac{1}{2},j}:=\frac{c^n_{p,i,j} + c^n_{p,i+1,j}}{2},\\
c^n_{p,i-\frac{1}{2},j}:=\frac{c^n_{p,i,j} + c^n_{p,i-1,j}}{2} .
\end{cases}
~
\text{and }\
\begin{cases}
\mathcal{A}^n_{p,i,j+\frac{1}{2}}:=\frac{\mathcal{A}^n_{p,i,j} + \mathcal{A}^n_{p,i,j+1}}{2},\\ 
\mathcal{A}^n_{p,i,j-\frac{1}{2}}:=\frac{\mathcal{A}^n_{p,i,j} + \mathcal{A}^n_{p,i,j-1}}{2},\\
\mathcal{A}^n_{p,i+\frac{1}{2},j}:=\frac{\mathcal{A}^n_{p,i,j} + \mathcal{A}^n_{p,i+1,j}}{2},\\
\mathcal{A}^n_{p,i-\frac{1}{2},j}:=\frac{\mathcal{A}^n_{p,i,j} + \mathcal{A}^n_{p,i-1,j}}{2} .
\end{cases}
\end{equation}
Further, the central differences are given by
\begin{equation}
\begin{cases}
{[c_{p,y}]}^n_{i,j+\frac{1}{2}} := \frac{c^n_{p,i,j+1} - c^n_{p,i,j}}{\Delta y},\\
{[c_{p,y}]}^n_{i,j-\frac{1}{2}} := \frac{c^n_{p,i,j} - c^n_{p,i,j-1}}{\Delta y},\\
{[c_{p,x}]}^n_{i+\frac{1}{2},j} := \frac{c^n_{p,i+1,j} - c^n_{p,i,j}}{\Delta x},\\
{[c_{p,x}]}^n_{i-\frac{1}{2},j} := \frac{c^n_{p,i,j} - c^n_{p,i-1,j}}{\Delta x}.
\end{cases}
\end{equation}
By using these notations the approximation for \re{$D_{p} \Delta c_{p}- \nabla \cdot [ c_{p} \mathcal{A}_{p}\left(x,t,\bu\left(\cdot,t\right)\right)]$} in \eqref{eq:TwoPop_Nonlocal_Chap5} is as follows: 
\begin{equation}
\begin{split}
\re{D_{p}\Delta c_{p}-\nabla \cdot [c_{p} \mathcal{A}_{p}\left(x,t,\bu\left(\cdot,t\right)\right)]} &= \\ &
 \simeq \frac{ \re{D_{p}} [c_{p,x}]^{n}_{i+\frac{1}{2},j} - \re{D_{p}} [c_{p,x}]^{n}_{i-\frac{1}{2},j} - c^{n}_{p,i+\frac{1}{2},j} \cdot \mathcal{A}^{n}_{p,i+\frac{1}{2},j} + c^{n}_{p,i-\frac{1}{2},j} \cdot \mathcal{A}^{n}_{p,i-\frac{1}{2},j}  }{ \Delta{x} }\\
& + \frac{ \re{D_{p}}[c_{p,y}]^{n}_{i,j+\frac{1}{2}} - \re{D_{p}} [c_{p,y}]^{n}_{i,j-\frac{1}{2}} - c^{n}_{p,i,j+\frac{1}{2}} \cdot \mathcal{A}^{n}_{p,i,j+\frac{1}{2}} + c^{n}_{p,i,j-\frac{1}{2}} \cdot \mathcal{A}^{n}_{p,i,j-\frac{1}{2}}  }{ \Delta{y} }.
\end{split}
\end{equation}
%
 \re{Next,} we shift our attention to the numerical approach of the adhesive flux $\mathcal{A}_{p}$ (that explores the effects of cell-cell and cell-matrix adhesion of cancer cells subpopulations). We carry out these computations off-grid by decomposing the sensing region $\textbf{B}\left(x,R\right)$ in \dt{$s\in \N^{*}$ annuli centred at $x$ (with the inner most central circle being of numerically negligible size), and for each annuli $k\in\{1,\dots,s\}$ (counted from the inner most to the outer most), we consider a radial decomposition of this in $2^{h+(k-1)}$ uniformly distributed radial sectors (with $h\in\N$ fixed), as illustrated in Figure \eqref{fig:Schematic_One} (for $h=2$). This leads to a collection of sectors $\{\mathcal{S}_{\nu}\}_{\nu=1,N_{s}}$, where 
 \[
 N_{s}:=\sum\limits_{k=1}^{s}2^{h+(k-1)}.
 \]}  
Then for each annulus sector $\mathcal{S}_{\nu}$, we evaluate the total primary cancer cell population $c_{1}$, the total mutated cancer cell population $c_{2}$ and the total ECM mass distributed on $\mathcal{S}_{\nu}$ that are given by 
\begin{equation}
\omega_{\mathcal{S}_{\nu},c_{1}} := \frac{1}{\lambda\left(\mathcal{S}_{\nu}\right)} \int\limits_{\mathcal{S}_{\nu}} c_{1}\left(x,t\right) \, dx, \quad \omega_{\mathcal{S}_{\nu},c_{2}} := \frac{1}{\lambda\left(\mathcal{S}_{\nu}\right)} \int\limits_{\mathcal{S}_{\nu}} c_{2}\left(x,t\right) \, dx,  \quad \text{and} \quad \omega_{\mathcal{S}_{\nu},v} := \frac{1}{\lambda\left(\mathcal{S}_{\nu}\right)} \int\limits_{\mathcal{S}_{\nu}} v\left(x,t\right) \, dx.
\end{equation}
\par 
\dt{Further}, on each sector $\mathcal{S}_{\nu}$, we \dt{consider} the off-grid barycentres of each annulus sector by $\textbf{b}_{\mathcal{S}_{\nu}}$ and \dt{the values of each of the three densities at these barycentres, namely $c_{1}(\textbf{b}_{\mathcal{S}_{\nu}},t)$ $c_{2}(\textbf{b}_{\mathcal{S}_{\nu}},t)$ and $v(\textbf{b}_{\mathcal{S}_{\nu}},t)$,  are obtained via interpolation with bilinear shape functions on the grid rectangles $\{ y^{i}_{\textbf{b}_{\mathcal{S}_{\nu}}} \}_{i=1,4}$ that contain $\textbf{b}_{\mathcal{S}_{\nu}}$. }Therefore, the approximation of the adhesion integral at each instance of time $t>0$ is given by 
\begin{equation}
\mathcal{A}_{p}\left(x,t,\bu\left(\cdot,t\right)\right) = \sum_{\nu=1}^{s2^{m}} \frac{\lambda \left(\mathcal{S}_{\nu}\right)}{R} \textbf{n}\left(\textbf{b}_{\mathcal{S}_{\nu}}\right) \cdot \mathcal{K}\left(\textbf{b}_{\mathcal{S}_{\nu}}\right) g_{\dt{p}}(\tilde{\bu}(\textbf{b}_{\mathcal{S}_{\nu}},t)), 
\end{equation}
where 
\[
\tilde{\bu}\left(\textbf{b}_{\mathcal{S}_{\nu}},t\right) := [\omega_{\mathcal{S}_{\nu},c_{1}}\left(t\right), \omega_{\mathcal{S}_{\nu},c_{2}}\left(t\right), \omega_{\mathcal{S}_{\nu},v}\left(t\right)]^{T},
\]
and \dt{$g_{p}(\tilde{\bu}(\textbf{b}_{\mathcal{S}_{\nu}},t))$ is the $p-$th component of}
\[
g\left(\tilde{\bu}\left(\textbf{b}_{\mathcal{S}_{\nu}},t\right)\right) = \left[ \textbf{S}_{cc} \left[\omega_{\mathcal{S}_{\nu},c_{1}}\left(t\right), \omega_{\mathcal{S}_{\nu},c_{2}}\left(t\right)\right]^{T} +  \dt{\omega_{\mathcal{S}_{\nu},v}(t)}\textbf{S}_{cv} \right]\cdot \left(K_{c} - \bc(\textbf{b}_{\mathcal{S}_{\nu}}, t) - v(\textbf{b}_{\mathcal{S}_{\nu}}, t)\right)^{+}.
\]
%
%
%
%
%
%
%
%
%
\subsection{Reconstruction of the Mutation Law in the Nonlocal Cancer Invasion Model}\label{Sect8}
Building on the approach described in local cancer invasion model Section \ref{Recon_Mutation_in_Local_Model}, we proceed now to address the reconstruction of the unknown mutation law $\Q(c_{1}, v)$ \dt{within the} nonlocal cancer invasion model \eqref{eq:TwoPop_Nonlocal_Chap5} from \dt{exact and noisy} measured data \dt{that are considered here again to be of the form prescribed in \eqref{comp_meas_1}-\eqref{noise_in_the_measurements}.}
%
%
\begin{figure}[!ht]
	\centering
        \setlength{\arrayrulewidth}{0.5pt}
        \begin{tabular}{lc}
        a)& \\
        &\hspace{-0.6cm} 
	\includegraphics[width=0.65\textwidth]{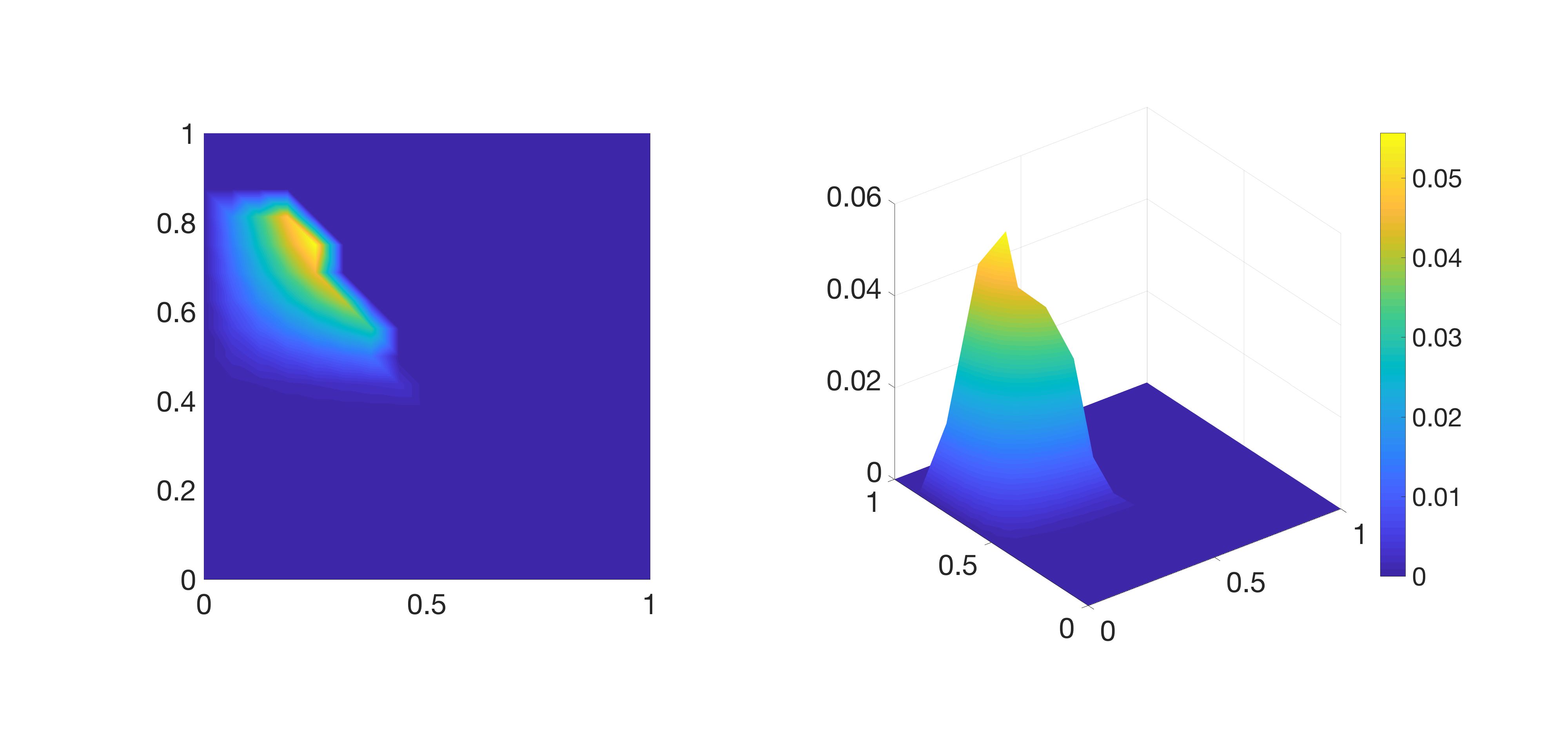}\\[-0.1cm]
         b)& \\
         &\hspace{-0.5cm}
	\includegraphics[width=0.3\textwidth]{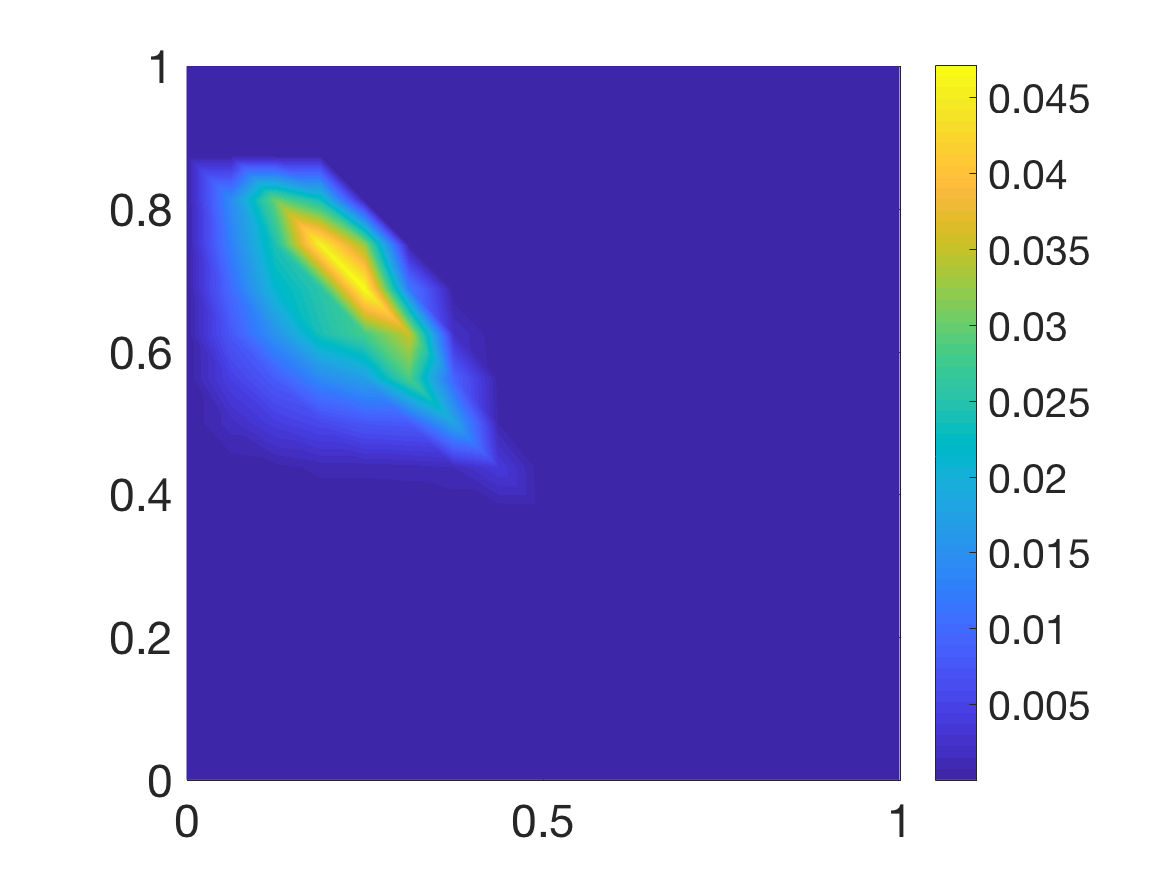}\hspace{-0.3cm}
	\includegraphics[width=0.3\textwidth]{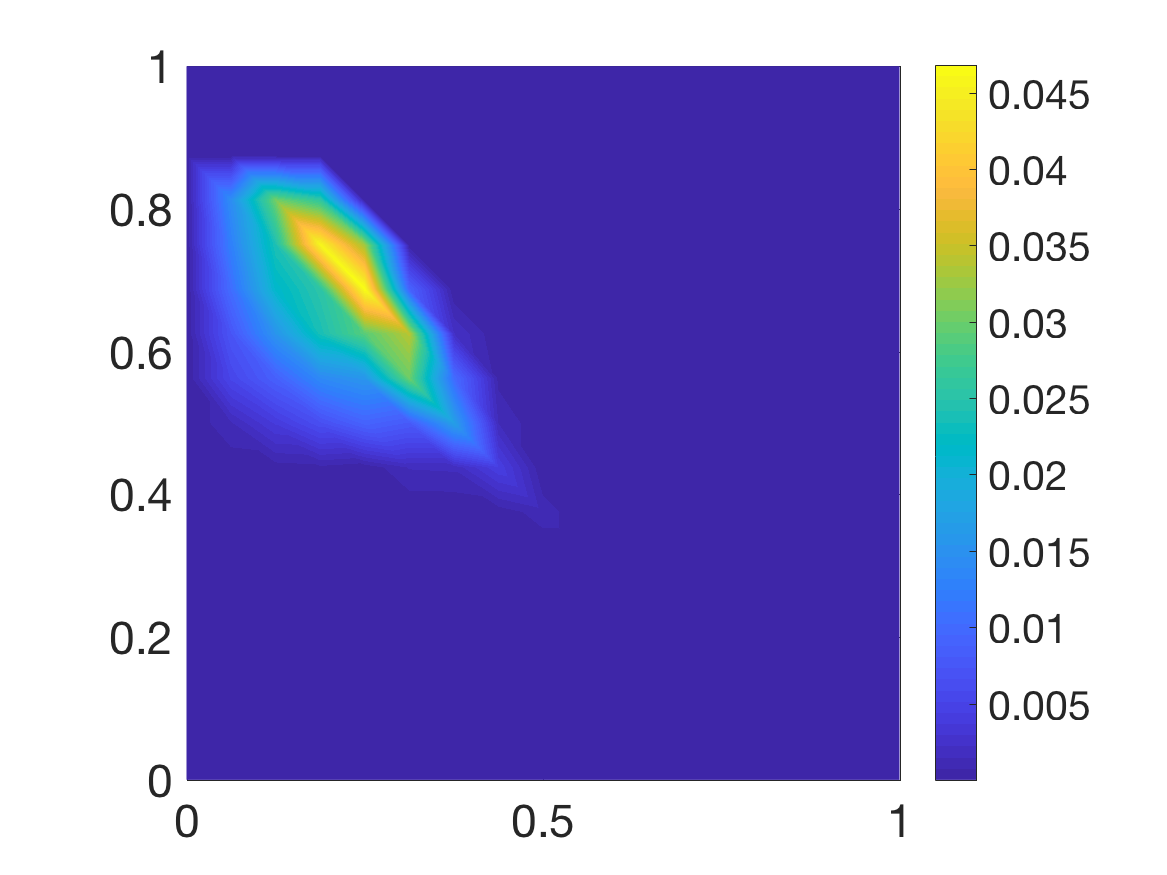}\hspace{-0.3cm}
	\includegraphics[width=0.3\textwidth]{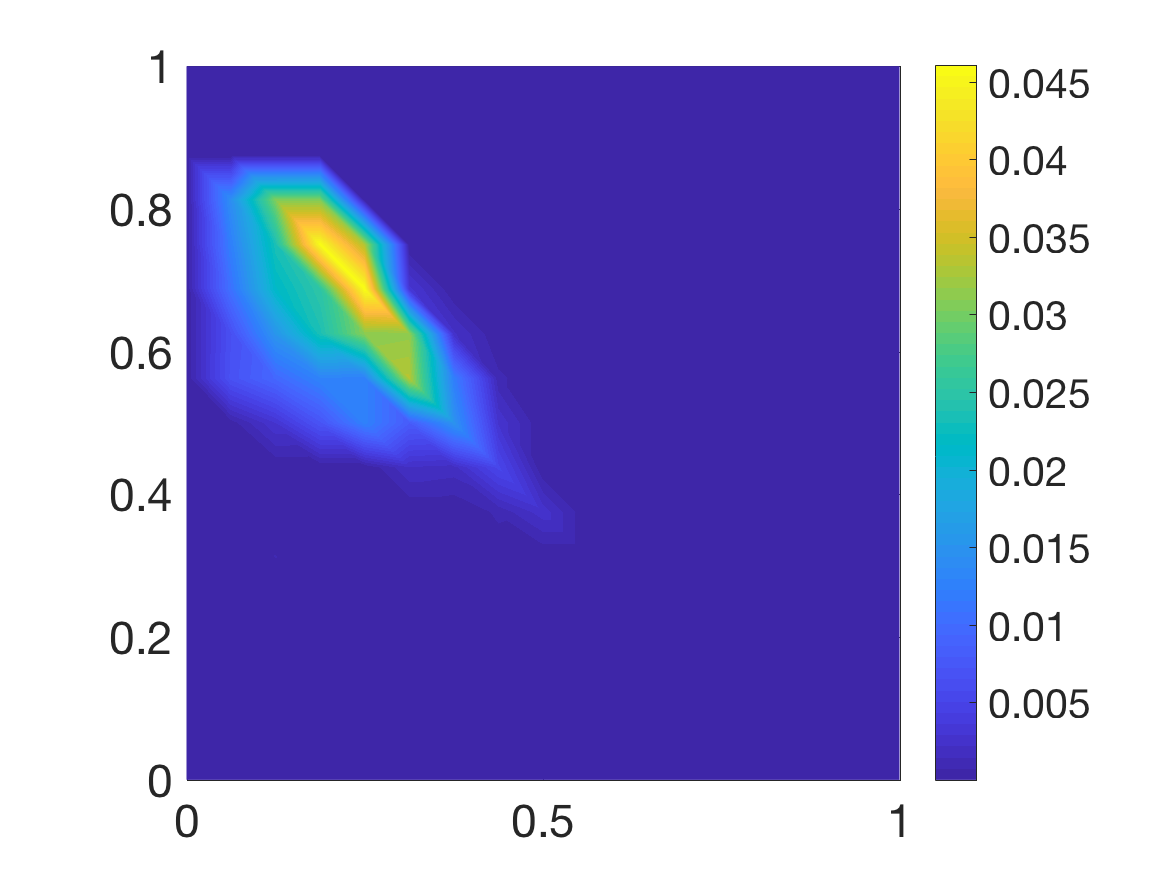}
	\end{tabular}
\caption{Reconstruction of \re{the general mutation law $Q(c_{1},v)$ for the nonlocal} model \eqref{eq:TwoPop_Nonlocal_Chap5}. Row a) the true mutation law restricted to $\A_{cv}$. Row b) the reconstructed  mutation law on $\A_{c\re{v}}$ in the presence of exact and noisy data: (left) exact data and $\alpha^{*}=10^{-4}$; (centre) $1\%$ noisy data and $\alpha^{*}=10^{-4}$; and (right) $3\%$ noisy data and $\alpha^{*}=10^{-4}$. For all plots in this figure: 1) first axis \re{shows} the values $c_{1}\in[\bar{c}_{1}^{min},\bar{c}_{1}^{max}]$; 2) second axis \re{shows} the values $v\in [\bar{v}^{min},\bar{v}^{max}]$; and 3) the colour bars represent the magnitude of mutation law or its reconstructions at each $(c_{1},v)\in \A_{cv}$. \re{The numerical simulations are obtained using the parameters given in Table \ref{paramSetTable}. }}
\label{fig:Reconstruction_of_Mutation_Law_NonLocal}
\end{figure}

Here, the exact \dt{measurement data that are given by
\bequ
 \label{exactMeasurements_nonlocal_model}
\tilde{c}^{*}_{1,exact}(x):= \bar{c}_{1}(x,t_{f}),\quad \tilde{c}^{*}_{2,exact}(x):= \bar{c}_{2}(x,t_{f}),\quad and \quad \tilde{v}^{*}_{exact}(x):=\bar{v}(x,t_{f}), \quad \forall x\in \Omega,
\eequ
where  $\bar{c}_{p}(x,t)$ $p=1,2$ and $\bar{v}(x,t)$ are the corresponding solution densities of the primary and mutated cancer cell populations as well as that of ECM that are obtained for model \eqref{eq:TwoPop_Nonlocal_Chap5} obtained when this uses as mutation law given by \eqref{eq:Mutation_DepOn_V} with the known term $\widetilde{\Q}_{2}(v)$ detailed in \eqref{explicit_measurement_case_ii}}. 

%
%
%
%
%
%
Figure \ref{fig:Reconstruction_of_Mutation_Law_NonLocal} shows the reconstruction of the cancer cell mutation law \re{$Q(c_{1},v)$} for model \eqref{eq:TwoPop_Nonlocal_Chap5} in the presence of  the measurements given by \eqref{exactMeasurements_nonlocal_model} that are considered here both exact and affected by a level of noise $\delta\in\{1\%, 3\%\}$. \re{As for the previous figures,} the first row shows the true mutation law restricted at the maximal accessible region  $\A_{cv}$ where the reconstruction is being attempted. The second row of the figure shows, from left to right, the reconstruction of the mutation law on $\A_{cv}$ with no noise, with $1\%$, and with $3\%$ of noise in the measured data, respectively. From this figure, we observe that we obtain a good mutation law reconstruction when the measurement data are not affected by noise. However, as expected, as soon as the level of noise increases in the measurements, the reconstruction gradually looses accuracy. 
%
%
%
%
\begin{figure}[!t]
\centering
\setlength{\arrayrulewidth}{0.5pt}
\begin{tabular}{lc}
a)&\\
& \hspace{-0.6cm}
\includegraphics[width=0.65\textwidth]{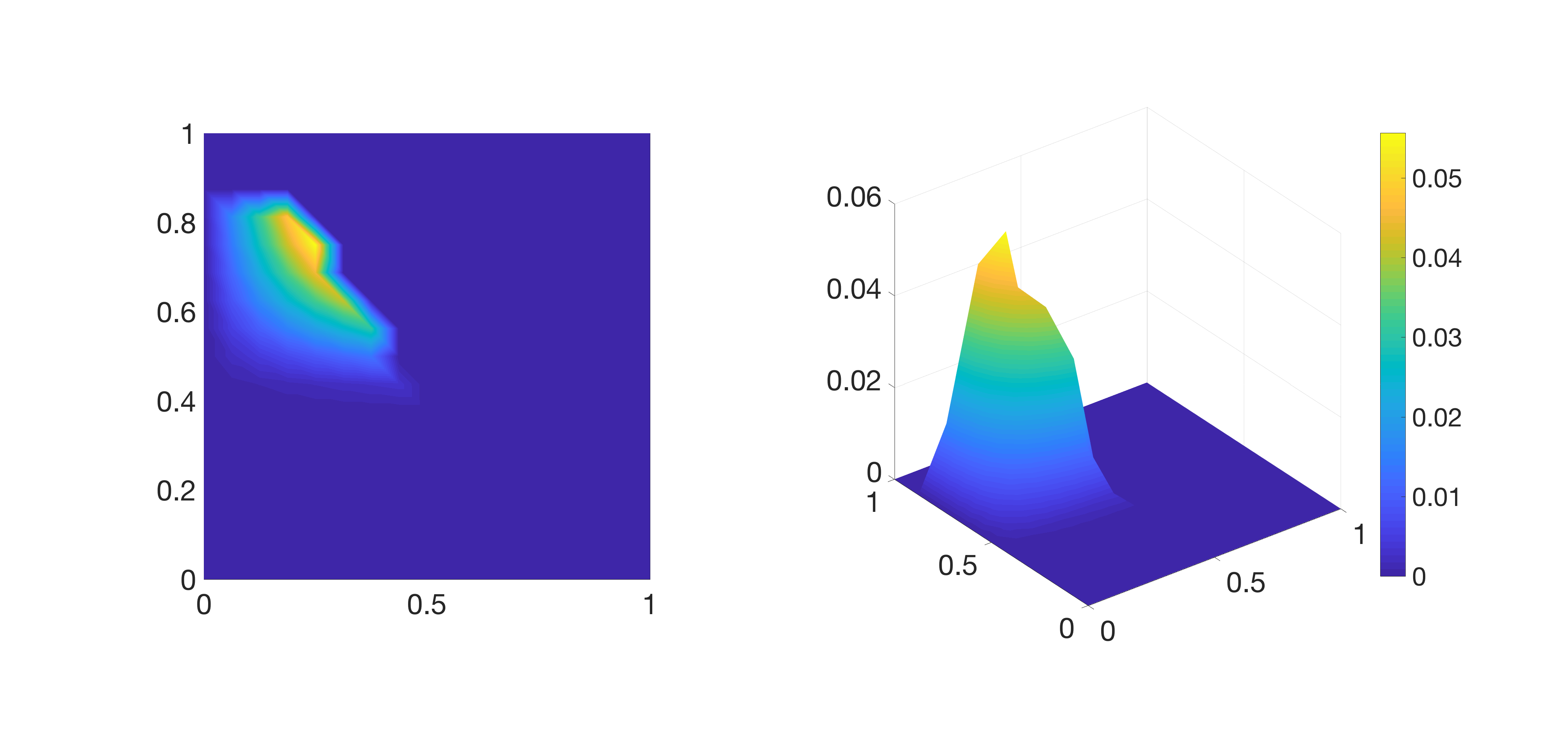}\\[-0.1cm]
b)&\\
& \hspace{-0.5cm}
\includegraphics[width=0.3\textwidth]{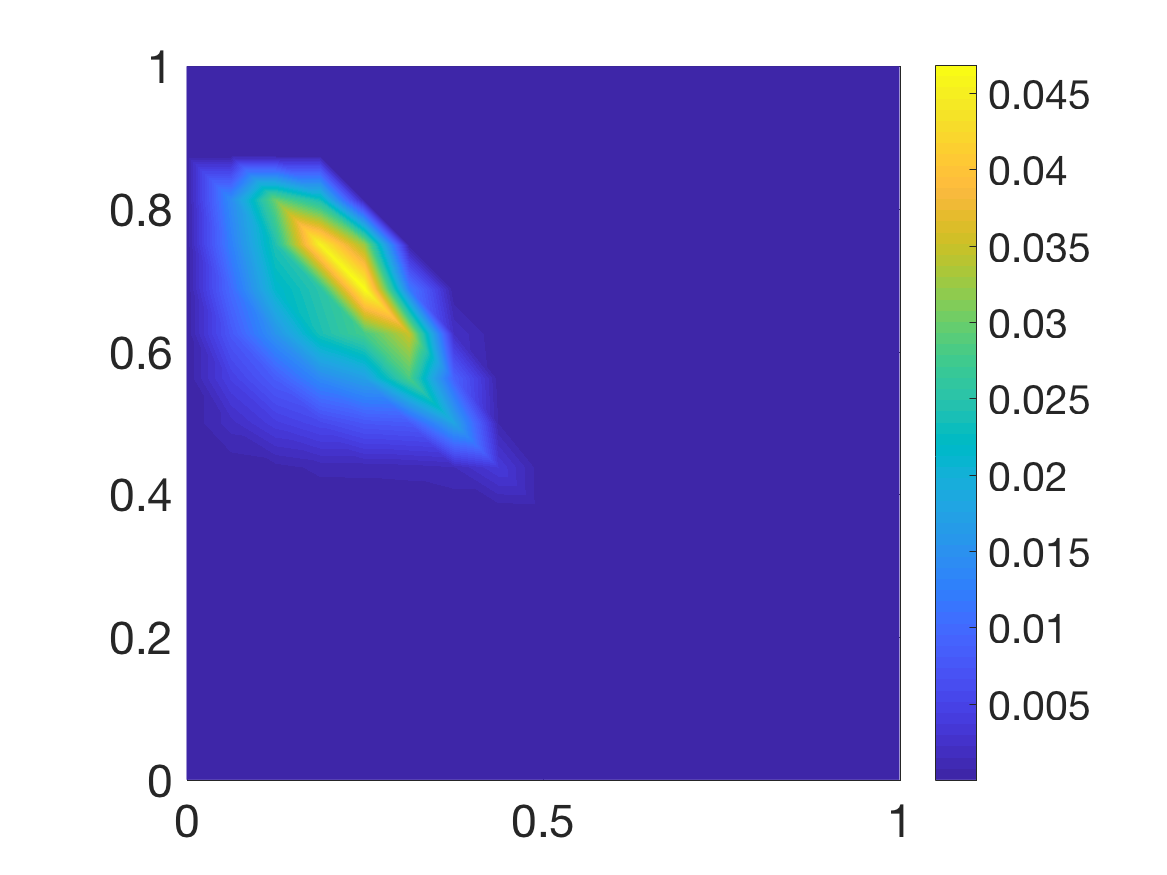}\hspace{-0.3cm}
\includegraphics[width=0.3\textwidth]{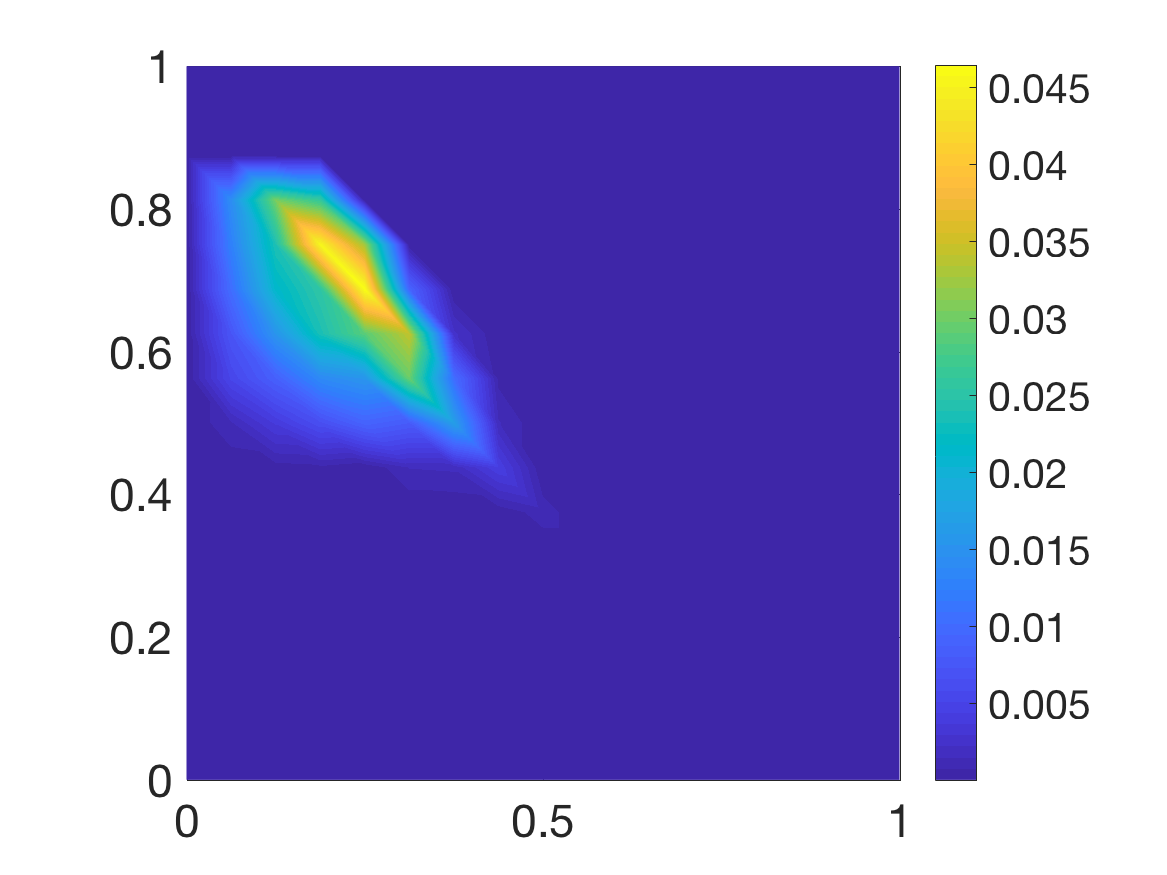}\hspace{-0.3cm}
\includegraphics[width=0.3\textwidth]{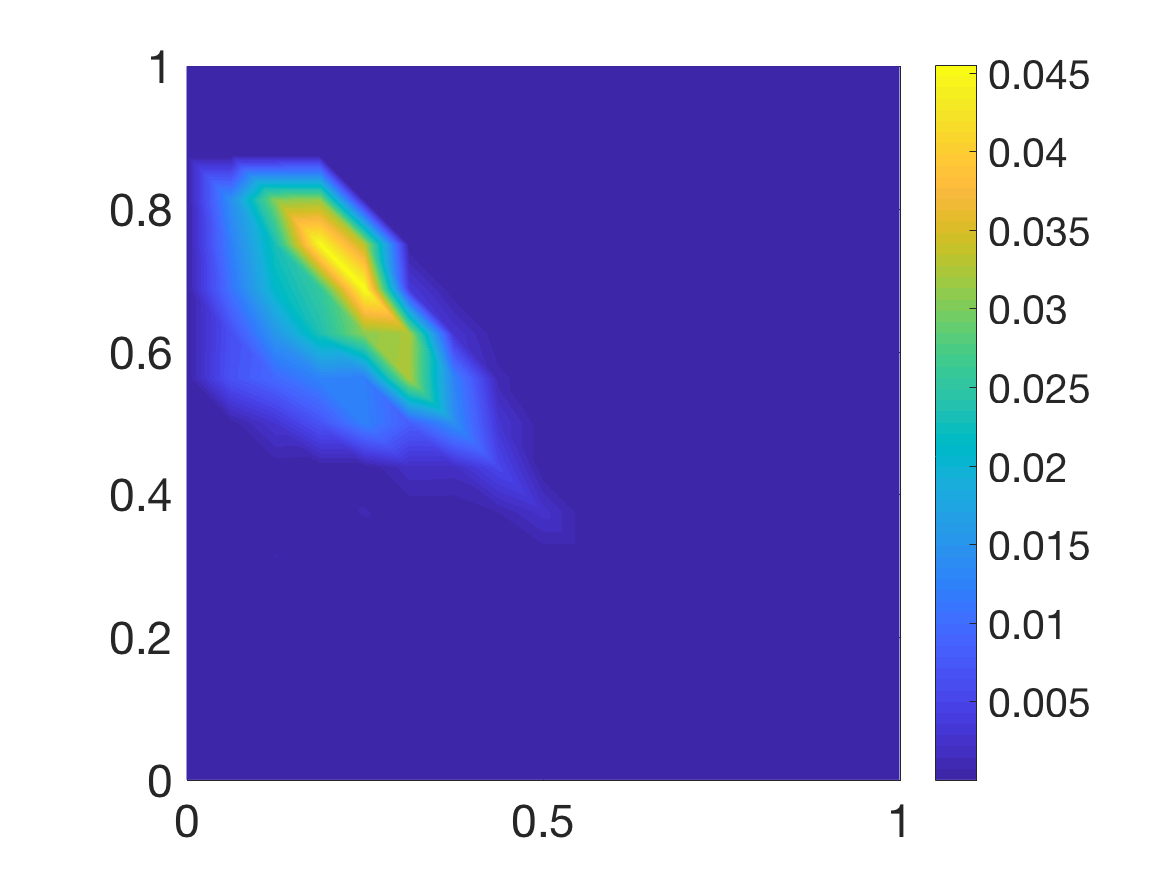}
\end{tabular}
\caption{Reconstruction of \re{the general mutation law $Q(c_{1},v)$} for the nonlocal model \eqref{eq:TwoPop_Nonlocal_Chap5} with mixed proliferation \re{rules for $c_{1}$ and $c_{2}$ cells}. Row a) the true mutation law restricted to $\A_{cv}$. Row b) the reconstructed mutation law on $\A_{cv}$ in the presence of exact and noisy data: (left) exact data and $\alpha^{*}=10^{-4}$; (centre) $1\%$ noisy data and $\alpha^{*}=10^{-4}$; and (right) $3\%$ noisy data and $\alpha^{*}=10^{-4}$. For all plots in this figure we have that: 1) the first axis represents the values for $c_{1}\in[\bar{c}_{1}^{min},\bar{c}_{1}^{max}]$; 2) second axis represents the values for $v\in [\bar{v}^{min},\bar{v}^{max}]$; and 3) the colour bars represent the magnitude of mutation law or its reconstructions at each $(c_{1},v)\in \A_{cv}$. \re{The parameters used for the numerical simulations are given in Table \ref{paramSetTable}. }}
\label{fig:Reconstruction_of_Mutation_Law_Mixed_Prol_NonLocal}
\end{figure}

\subsection{Reconstruction of the Mutation Law \re{for a Different Cell Proliferation Rules for $c_{2}$ Population}}

\re{As in Section~\ref{SectLocalMixedGrowth}, here we investigate the reconstruction of the general mutation law $Q(c_{1},v)$ when cancer cell growth is described by different rules: logistic proliferation for the $c_{1}$ cells, and Gompertz proliferation for the $c_{2}$ cells. Again, we do not see any significant differences between the reconstruction of mutation law in this case (see Figure \ref{fig:Reconstruction_of_Mutation_Law_Mixed_Prol_NonLocal}) and in the previous case where both cell populations grow logistically (see Figure \ref{fig:Reconstruction_of_Mutation_Law_NonLocal}). We believe that is because in both cases, the population that starts mutating (i.e., the $c_{1}$ population) grows logistically. }

%
%

%
\section{Conclusion}\label{Sect:Conclusion}
In this work we explored a new inverse problem that addresses the reconstruction of the cancer cells mutation law in \re{two heterotypic cancer invasion models: a model with a local cancer cell haptotaxis flux towards ECM (see equations \eqref{eq:Direct_Mutation_DepOn_C}), and a model with a nonlocal haptotaxis flux generated by cell-cell and cell-ECM adhesion forces (see equations \eqref{eq:TwoPop_Nonlocal_Chap5}). For the reconstruction of this mutation law through an inverse problem Tikhonov regularisation-based approach, we used a numerically-generated spatial tumour snapshot data assumed to be acquired at a later stage in the tumour evolution (in practice, the data can be provided through a medical imaging scan).}  

This inverse problem approach \re{was} implemented computationally via a mixed finite differences - finite element numerical scheme. Specifically, on one hand, we used a Crank-Nicholson-type finite difference scheme for the  discretisation of the \emph{forward models} that arises in each of the two tumour invasion dynamics considered here (i.e., local cancer cell invasion, and nonlocal cancer cell invasion). On the other hand, we develop\re{ed} a finite element approach involving a bilinear shape functions on a square mesh for the discretisation of mutation law candidates recruited from a proposed space of functions $\S$, as well as \re{for} their evaluation on a maximal accessible regions where the mutation law reconstruction \re{was} performed. Finally, these two parts \re{were} appropriately assembled in an optimisation solver that \re{sought} to reconstruct the cancer cell mutation law by minimising over the space $\S$ the emerging Tikhonov functionals that \re{were }formulated in each of the two cases considered.

Finally, this inversion approach was explored and tested on the reconstruction \re{of different} cancer cell mutation laws used in local and nonlocal cancer modelling: (i) reconstructing the cancer cell mutation \re{assuming a linear dependence only on the $c_{1}$ cell population;} (ii) reconstructing the cancer cell mutation \re{assuming a linear dependence on the $c_{1}$ population and a nonlinear dependence on the ECM density}, and (iii) reconstructing \re{a very general mutation law, while assuming no prior knowledge about the mutation. All numerical reconstruction results for the local and nonlocal models showed that} for exact measurements we obtain\re{ed} a good reconstruction of the mutation law, while for increasingly noisy measurements the reconstruction gradually deteriorates.
\appendix
\section{Parameters used in computations}\label{modelParams}
For all the cancer cells mutation laws reconstructions considered in this work, we use the parameter set specified in Table \ref{paramSetTable}.
\begin{table}[H]
	\centering
	\caption{Summary of parameter values that have used for two local and nonlocal sub-population of cancer cells.}
	\vspace{0.5\baselineskip}
	\begin{tabular}{cclc}
		\hline
		Parameter                    &Value                                  & Description                                                         & Reference \\
		\hline
		$D_{1}$                        &$0.00675$                         & diffusion of primary tumour                                 &\cite{Domschke_et_al_2014}  \\
		
    	$D_{2}$                        &$0.00675$                         & diffusion of secondary tumour                            &\cite{Domschke_et_al_2014}  \\
		
		$\eta_{1}$                     &$2.85\times 10^{-2}$         & haptotaxis to ECM from $c_{1}$                           &\cite{Peng_et_al_2017}  \\
		
    	$\eta_{2}$                     &$2.85\times 10^{-2}$         & haptotaxis to ECM from $c_{2}$                           &\cite{Peng_et_al_2017}  \\
		
		$\mu_{c}$               & $0.25$                                     & proliferation of tumour cells c                            &  \cite{Shuttleworth_2018_two} \\
		
                $K_{c}$              & $1$                                     & tissue carrying capacity                                     & \cite{Dumitru_et_al_2013}  \\
		
		$\rho$                 & $2$                                          & ECM degradation coefficient                                      &\cite{Shuttleworth_2019}  \\
		
		$\mu_{v}$              & $0.40$                                     & ECM remodelling coefficient                                        & \cite{Dumitru_et_al_2013}  \\
		
		\dt{$\textbf{S}_{cc}$}  & \dt{$\begin{pmatrix} 0.5 & 0 \\ 0   & 0.3 \end{pmatrix}$} & \dt{cell-cell adhesion function}  & \dt{\cite{Domschke_et_al_2014}}\\[0.4cm]
		
		\dt{$\textbf{S}_{cv}$}  & \dt{$\begin{pmatrix} 0.1 & 0 \\ 0   & 0.5 \end{pmatrix}$} & \dt{cell matrix adhesion function}  & \dt{\cite{Domschke_et_al_2014}}\\
		
		$t_{1,2}$          & $10$                                        & time initiation for mutations                           &  \cite{Shuttleworth_2018_two,Andasari_2011}  \\	
		
		\dt{$t_{s}$}	&\dt{$3$}                                                 & \dt{time-steepness coefficient} 				&   \dt{\cite{Shuttleworth_2018_two}} \\
		
		\dt{$\delta_{0}$}          & $0.3$                                        & mutation from primary tumour                            &  \cite{Shuttleworth_2018_two,Andasari_2011}  \\
		
		$\Delta{x}$             &$0.03125$                               & discretisation step size for $\G_{_{\Omega}}$                 & \cite{Dumitru_et_al_2013}  \\
		
		$\Delta t$             &$ 10^{-3}$                               & time step size                  & \cite{Dumitru_et_al_2013}  \\

		$\Delta{\eta}$        &$ 0.0625$                                & mesh size used for \dt{$\G^{1}_{_{M}}$   and $\G^{2}_{_{M}}$}                                                              &     Estimated        \\
				
		\hline
	\end{tabular}
	\label{paramSetTable}
\end{table}
%
%
%
%
%
%
%
%
%
%
%
\section*{Acknowledgement}
The first author would like to thank Saudi Arabian Cultural Bureau in the United Kingdom (UKSACB) on behalf of Taif University and the University College of Al-Khurmah in Saudi Arabia for supporting and sponsoring his PhD studies at the University of Dundee.\\\\
\bibliographystyle{aims}
\bibliography{paper_Two_bibiliography}
%
%
%
%
%
%
%
\end{document}